\documentclass[a4paper,reqno]{amsart}

\pdfoutput=1

\usepackage{nicefrac,amsmath,amssymb,amsthm,url,verbatim,enumitem,stmaryrd,mathtools,microtype,bbm,scalerel}
\usepackage[utf8]{inputenc}
\usepackage[T1]{fontenc}
\usepackage{libertine}
\usepackage{dsfont}
\usepackage[libertine,cmintegrals,cmbraces]{newtxmath}
\usepackage[dvipsnames,svgnames,x11names]{xcolor}
\definecolor{darkred}{RGB}{139,0,0}
\definecolor{darkblue}{RGB}{0,0,139}
\definecolor{darkgreen}{RGB}{0,100,0}
\usepackage[mathscr]{euscript}
\usepackage[pagebackref]{hyperref}
\hypersetup{
	colorlinks   = true, 
	urlcolor     = darkred, 
	linkcolor    = darkred, 
	citecolor   = darkgreen}
\usepackage{tikz}
\usepackage{tikz-cd}
\usetikzlibrary{calc}

\tikzset{
  snake arrow/.style={
    rounded corners,
    to path={ -- ([xshift=2ex]\tikztostart.east)
              |- ($(\tikztostart.south)!.5!(\tikztotarget.north)$) \tikztonodes
              -| ([xshift=-2ex]\tikztotarget.west)
              -- (\tikztotarget)}
  }
}
\usetikzlibrary{decorations.markings,patterns,arrows.meta}
\tikzset{>={Straight Barb[scale=0.8]}}
\tikzcdset{
	arrow style=tikz,
	diagrams={>={Straight Barb[scale=0.8]}}
}
\usepackage[capitalise,noabbrev]{cleveref}

\calclayout

\AtBeginDocument{%
   \def\MR#1{}
}

\setenumerate[1]{label=(\roman*),}
\setitemize[1]{label=\raisebox{0.25ex}{\tiny$\bullet$}}

\setlist[enumerate]{leftmargin=*}
\setlist[itemize]{leftmargin=*}

\newcommand\smallsquare{{\mathbin{\text{\raise0.17ex\hbox{\scalebox{.7}{$\blacksquare$}}}}}}

\newcommand{\circled}[1]{\raisebox{.5pt}{\textcircled{\raisebox{-.9pt} {#1}}}}

\newcommand{\BO}{\ensuremath{\mathrm{BO}}}

\newcommand{\Diff}{{\ensuremath{\mathrm{Diff}}}}

\newcommand{\BDiff}{\ensuremath{\mathrm{BDiff}}}

\newcommand{\Map}{\ensuremath{\mathrm{Map}}}

\newcommand{\Emb}{\ensuremath{\mathrm{Emb}}}

\newcommand{\BEmb}{\ensuremath{\mathrm{BEmb}}}

\newcommand{\interior}{\ensuremath{\mathrm{int}}}

\newcommand{\inc}{\ensuremath{\mathrm{inc}}}

\newcommand{\res}{\ensuremath{\mathrm{res}}}

\newcommand{\Fun}{\ensuremath{\mathrm{Fun}}}

\newcommand{\Ind}{\ensuremath{\mathrm{Ind}}}

\newcommand{\PSh}{\mathrm{PSh}}

\DeclareMathAlphabet{\mathpzc}{OT1}{pzc}{m}{it}
\newcommand{\catsingle}[1]{\ensuremath{\mathscr{#1}}}

\newcommand{\icat}[1]{\ensuremath{\mathscr{#1}}}

\newcommand{\Cosp}{\ensuremath{\mathrm{Cosp}}}

\DeclareMathOperator*{\colim}{\ensuremath{\mathrm{colim}}}

\newcommand{\hocolim}{\ensuremath{\mathrm{hocolim}}}

\newcommand{\holim}{\ensuremath{\mathrm{holim}}}

\newcommand{\id}{\mathrm{id}}

\newcommand{\fib}{\mathrm{fib}}

\newcommand{\op}{\mathrm{op}}

\newcommand{\Pro}{\mathrm{Pro}}

\newcommand{\forget}{\mathrm{forget}}

\newcommand{\comp}{{\mathrm{comp}}}

\newcommand{\ManInf}{\ensuremath{\icat{M}\mathrm{fld}}}

\newcommand{\BordInf}{\ensuremath{\icat{B}\mathrm{ord}}}

\newcommand{\ncBordInf}{\ensuremath{\mathrm{nc}\icat{B}\mathrm{ord}}}

\newcommand{\Cat}{\ensuremath{\icat{C}\mathrm{at}}}

\newcommand{\CatInf}{\ensuremath{\icat{C}\mathrm{at}_\infty}}

\newcommand{\COSPAN}{\ensuremath{\mathrm{COSPAN}}}

\newcommand{\tw}{\ensuremath{\mathrm{tw}}}
\newcommand{\fr}{\ensuremath{\mathrm{fr}}}

\newcommand{\Exh}{\mathrm{Exh}}
\newcommand{\open}{\mathrm{open}}
\newcommand{\cpct}{\mathrm{cpct}}

\newcommand{\OpenMan}{\ensuremath{\icat{M}\mathrm{fld}}}
\newcommand{\CpctMan}{\ensuremath{\icat{M}\mathrm{fld}^{\mathrm{cpt}}}}

\newcommand{\ncMan}{\mathrm{nc}\ensuremath{\icat{M}\mathrm{fld}}}

\newcommand{\bordEnd}{\mathrm{End}_{\mathrm{bord}}}
\newcommand{\proEnd}{\mathrm{End}_{\mathrm{pro}}}

\newcommand{\fix}{\mathrm{fix}}
\newcommand{\free}{\mathrm{free}}

\newtheorem{bigthm}{Theorem}

\newtheorem{thm}{Theorem}[section]
\newtheorem*{nthm}{Theorem}
\newtheorem{lem}[thm]{Lemma}

\newtheorem{cor}[thm]{Corollary}

\theoremstyle{definition}
\newtheorem{dfn}[thm]{Definition}

\theoremstyle{remark}

\newtheorem{rem}[thm]{Remark}

\newtheorem{construction}[thm]{Construction}

\begin{document}

\title{A bordism perspective on ends of manifolds}

\author{Fadi Mezher}
\address{Department of Mathematics, Karlsruhe Institute of Technology, 76131 Karlsruhe, Germany}
\email{fadi.mezher@kit.edu}

\begin{abstract}
We refine the classical pro-homotopy type of the end of a manifold by recording the boundaries of a compact exhaustion together with the bordisms between them. We show that this assignment, the \emph{bordism-end}, is functorial in embeddings, and independent of the choice of exhaustion. Its target is the free sequential colimit completion of the bordism $\infty$-category. Along the way of defining the bordism-end, we show that the category of all smooth $d$-dimensional manifolds and embeddings is formally obtained from the analogous category of compact manifolds. Through the formalism of the bordism-end, it follows that a manifold is diffeomorphic to the interior of a compact manifold if and only if its bordism-end is representable. Moreover, this formalism yields a fibre sequence relating the classifying spaces of the compactly supported diffeomorphism group of a manifold and its full diffeomorphism group, and a space of bordism-ends. Finally, we set up the pro-end functorially in embeddings, with target the cospan $\infty$-category of pro-spaces, and show that it factors through the bordism-end.
\end{abstract}

\maketitle

\section{Introduction}
Given a smooth manifold $M$ with empty boundary (potentially non-compact), a choice of a proper Morse function on $M$ yields a compact exhaustion of $M$, that is, a sequence of codimension-0 smooth compact submanifolds $\{K_i\}_{i\in\mathbb{N}}$, such that $K_i \subset \interior(K_{i+1})$, and such that $M =\bigcup_i K_i$. Classically (see for instance \cite{SiebenmannThesis}), one records the "ends" of $M$ as the inverse system $\{M\setminus K_i\}$, as an object in the pro-category of spaces. In this work, we suggest a bordism theoretic refinement of the classical treatment of ends, motivated by the following:  we can consider the sequential system $\{\partial K_i\}_i$ as an $\Ind$-object in the bordism category, where the morphisms $\partial K_i \rightsquigarrow \partial K_{i+1}$ are the bordisms $K_{i+1} \setminus\interior(K_i)$. More precisely, the following represents the list of objectives we work towards:
\begin{enumerate}[label=Objective (\roman*)]
    \item Assemble the assignments $M \mapsto\{M\setminus K_i\}$, and $M\mapsto \{\partial K_i\}$ into functors from the category of smooth $d$-manifolds and smooth embeddings between them, into convenient categories where $\{M\setminus K_i\}$ and $\{\partial K_i\}$ naturally live. We would also like to relate these two functors. As an intermediate step, we;
    \item Show that the category of all smooth $d$-manifolds (with empty boundary) can be formally obtained from the category of smooth, compact $d$-manifolds (with boundary);
    \item Characterize when a manifold is diffeomorphic to the interior of a compact manifold, as a condition on its bordism-end;
    \item Relate the classifying spaces of the compactly supported diffeomorphism group of a manifold and the full diffeomorphism group, and spaces of bordism-ends by means of a fibre sequence.
\end{enumerate}

We begin by stating our answer to the first objective above. For $\catsingle{C}$ an $\infty$-category, we let $\PSh^\mathbb{N}(\catsingle{C})$ denote the \emph{free sequential colimit completion} of $\catsingle{C}$: this plays the same role for sequential colimits as $\PSh(\catsingle{C})$ plays for all colimits. More concretely, $\PSh^\mathbb{N}(\catsingle{C})$ contains the full subcategory of $\PSh(\catsingle{C})$ on sequential colimits of representables (and in fact is equal to it, see Lemma \ref{lem:naive-colim}). Thus, any sequential diagram $c_1 \to c_2 \to c_3 \cdots$ in $\catsingle{C}$ gives rise to an object in $\PSh^\mathbb{N}(\catsingle{C})$: for $\catsingle{C}$ the bordism $\infty$-category, the former will be the category where the bordism-end of a manifold naturally lives (we refer to §\ref{sec:definition} for more details on sequential colimit completions). Furthermore, for $\catsingle{C}$ a category admitting pushouts, let $\Cosp(\catsingle{C})$ denote the cospan $\infty$-category informally described as the category with objects the objects of $\catsingle{C}$, and where a morphism $c \rightsquigarrow d$ is the datum of a cospan $c\rightarrow h \leftarrow d$, and where composition is obtained by pushouts. Finally, we let $\OpenMan_d$ denote the $\infty$-category of smooth, $d$-manifolds with empty boundary and spaces of smooth embeddings between them, and $\BordInf_d$ denote the smooth, $d$-dimensional bordism $\infty$-category. We can now state our answer to the first item of the above list (in the text as Theorem \ref{thm:bordEnd-proEnd}). 
\begin{bigthm}\label{thmA}
    There are two functors
    \begin{center}
        $\bordEnd \colon \OpenMan_d \to \PSh^\mathbb{N}(\BordInf_d)$,\quad $\proEnd\colon \OpenMan_d \to \Cosp(\Pro(\catsingle{S}))$
    \end{center}
    and an arrow making the following triangle commute
    \begin{center}
        \begin{tikzcd}
            & \PSh^\mathbb{N}(\BordInf_d) \arrow[dashed]{dd}
            \\
            \OpenMan_d \arrow{ur}{\bordEnd} \arrow{dr}[swap]{\proEnd}
            &
            \\&\Cosp(\Pro(\catsingle{S}))
        \end{tikzcd}
    \end{center}
    Furthermore, for a manifold $W$ and a choice $\{K_i\}$ of a compact exhaustion, $\bordEnd(W) \simeq \{\partial K_i\}$, and $\proEnd(W) \simeq \{W\setminus K_i\}$.
\end{bigthm}
\subsection{From compact to all manifolds}
We begin by describing our answer to Objective (ii). This step is of independent interest, and proves crucial in our answers to the remaining objectives. Let $\CpctMan_d$ denote the $\infty$-category of compact manifolds, potentially with boundary, and spaces of smooth embeddings between them (where we impose no boundary condition). There is a functor $\interior\colon \CpctMan_d \to \OpenMan_d$, sending a compact manifold to its interior, and as any manifold is isotopy equivalent to its interior, $\interior$ is a fully faithful functor. Finally, the category $\ManInf_d$ admits sequential colimits: indeed, for a sequence $M_0 \hookrightarrow M_1 \hookrightarrow\cdots$ of smooth manifolds with empty boundary and smooth embeddings between them, one can form the colimit of that system in topological spaces, check that the smooth embeddings endow this colimit with the structure of a smooth manifold, and check the universal property of sequential colimits (which is equivalent to the cofunctor $\Emb(-,N)$ being good in the sense of Weiss \cite[1.4]{Weiss}).

Our answer to Objective (ii) takes the following form (in the text as Theorem \ref{thm:openMfldFreeColim}):
\begin{bigthm}\label{thmB}
    The unique, sequential colimit preserving extension $\interior_!$ of the functor $\interior$ defines an equivalence of $\infty$-categories
    $$\interior_! \colon \PSh^\mathbb{N}(\CpctMan_d) \xrightarrow[]{\simeq} \OpenMan_d$$
    Its inverse is obtained by a restricted Yoneda functor
    $$\Exh\colon \OpenMan_d \hookrightarrow\PSh^\mathbb{N}(\ManInf_d) \xrightarrow[]{\interior^\ast} \PSh^\mathbb{N}(\CpctMan_d)$$
\end{bigthm}
We emphasize that the pullback functor $\interior^\ast$ does not exist formally; one is required to check that the pullback functor $\interior^\ast \colon \PSh(\ManInf_d) \to \PSh(\CpctMan_d)$ restricts to a functor on $\PSh^\mathbb{N}$, which is carried out in Lemma \ref{lem:restriction-functor-exists}. The above theorem achieves two things simultaneously:
\begin{itemize}
    \item It assembles the association $M \mapsto \{K_i\}$, sending a manifold to its compact exhaustion, into a functor; indeed, $\Exh(M)$ will be equivalent to $\{K_i\}$ in $\PSh^\mathbb{N}(\CpctMan_d)$.
    \item It gives a formal construction of all manifolds with empty boundary (including open manifolds) from compact manifolds. It also satisfies the following universal property: given a functor $F\colon \CpctMan_d \to \catsingle{C}$ where $\catsingle{C}$ admits sequential colimits, there is a unique, sequential colimit preserving extension $F_! \colon \ManInf_d \to \catsingle{C}$.
\end{itemize}
\subsection{Bordism-ends}
The main object of this work is the functor $\bordEnd\colon\OpenMan_d \to \PSh^\mathbb{N}(\BordInf_d)$; indeed, all of our results center around that functor, and the functor $\proEnd$ will be defined through it. We now describe our definition of this functor, and highlight how Objectives (iii) and (iv) follow (almost) formally from it.

\subsubsection{Objective (i), part 1: defining $\bordEnd$}
We give a quick outline of how to define $\bordEnd$, which is one of the aims we highlighted in Objective (i). A crucial input in this definition is an equivalence of categories relating the category $\CpctMan_d$ to a slice of the bordism $\infty$-category: by \cite[Corollary A.2]{AdelaJan}, there is an equivalence of categories
$$\Phi\colon\CpctMan_d \simeq \left(\BordInf_d\right)_{\emptyset/}$$
which informally views a compact $d$-manifold $M$ as a bordism $\emptyset \rightsquigarrow\partial M$ (where $\partial M$ may be empty). From Theorem \ref{thmB}, it follows that there is an equivalence
$$\PSh^\mathbb{N}(\Phi) \colon \ManInf_d \simeq \PSh^\mathbb{N}\left((\BordInf_d)_{\emptyset/}\right)$$
It can be shown that $\PSh^\mathbb{N}$ behaves well with undercategories (Lemma \ref{lem:slice}), and therefore we obtain an equivalence
$$\ManInf_d \simeq \PSh^\mathbb{N}(\BordInf_d)_{y(\emptyset)/}$$
where $y$ denotes the Yoneda embedding. We will define the boundary, and bordism-end functors (in the text as Definition \ref{def:boundary-bordismends-functors}), as

\begin{itemize}
    \item $\partial\colon \CpctMan_d \to \BordInf_d$, as the composite $\CpctMan_d \simeq (\BordInf_d)_{\emptyset/} \xrightarrow[]{\mathrm{fgt}}\BordInf_d$;
    \item $\bordEnd\colon \ManInf_d\to \PSh^\mathbb{N}(\BordInf_d)$, as the composite 
    $$\ManInf_d\simeq \PSh^\mathbb{N}(\BordInf_d)_{y(\emptyset)/}\xrightarrow[]{\mathrm{fgt}}\PSh^\mathbb{N}(\BordInf_d)$$
\end{itemize}

\subsubsection{Objective (iii): a pullback square} One formal consequence of the above definition is Objective (iii). Indeed, we obtain the following pullback square (in the text as Theorem \ref{thm:Ind-End}).
\begin{bigthm}\label{thmC}
    For all $d\in\mathbb{N}$, there is a pullback square of $\infty$-categories
    \begin{center}
        \begin{tikzcd}
            \CpctMan_d \arrow[hook]{r}{\interior} \arrow{d}[swap]{\partial} & \ManInf_d\arrow{d}{\bordEnd}
            \\
            \BordInf_d \arrow[hook]{r}[swap]{y} & \PSh^\mathbb{N}(\BordInf_d)
        \end{tikzcd}
    \end{center}
\end{bigthm}
In more concrete terms, this square is saying the following. Given a manifold $M$, one can consider its bordism-end $\bordEnd(M)$; if that object is equivalent to a representable in $\PSh^\mathbb{N}(\BordInf_d)$, say $y(P)$ for $P$ a closed $(d-1)$-dimensional manifold, then $M$ is isotopy equivalent to the interior of a compact manifold with boundary $P$, and $P$ is only well defined up to invertible bordisms. Finally, one can show that isotopy invertible embeddings are equivalent to diffeomorphisms (Lemma \ref{lem:invertible-emb-diff}), which in particular allows us to replace "isotopy equivalent" by "diffeomorphic" in the previous sentence. This pullback square thus satisfies our objective (iii).
\subsubsection{Objective (iv): a fibre sequence}
The formalism of the pullback square of Theorem \ref{thmC} allows us now to show the following, as an answer to objective (iv) (in the text as Theorem \ref{thm:compactly-supported-diff}).
\begin{bigthm}\label{thmD}
    Let $Z \in \PSh^\mathbb{N}(\BordInf_d)$. Then, there is a pullback square of spaces
    \begin{center}
        \begin{tikzcd}
            \bigsqcup_{[W,\bordEnd(W)\simeq Z]}\BDiff_c(W) \arrow{r} \arrow{d} & \bigsqcup_{[W]}\BDiff(W)\simeq (\ManInf_d)^\simeq \arrow{d}{\bordEnd^\simeq}
            \\
            \ast \arrow{r}[swap]{\{Z\}} & \PSh^\mathbb{N}(\BordInf_d)^\simeq
        \end{tikzcd}
    \end{center}
    where
    \begin{itemize}
        \item $[W,\bordEnd(W)\simeq Z]$ denotes the set of pairs of a manifold and an identification $\bordEnd(W)\simeq Z$, up to diffeomorphisms rel. the bordism-end identification;
        \item $\Diff_c(W)$ denotes the space of compactly supported diffeomorphisms of such $W$;
        \item $[W]$ denotes the set of diffeomorphism classes of manifolds in $\ManInf_d$.
        \item $(-)^\simeq$ denotes the groupoid core: informally, if $\catsingle{C}$ is an $\infty$-category, then $\catsingle{C}^\simeq$ is the $\infty$-groupoid (alias space) obtained by discarding all non-invertible morphisms in $\catsingle{C}$.
    \end{itemize}
\end{bigthm}
Applying based loop spaces to the above fibre sequence, we obtain a fibre sequence relating $\Diff_c(W)$ to $\Diff(W)$; this is reminiscent of the fibre sequence set up by Cardim \cite[Theorem A]{Cardim}. However, in \emph{loc. cit.}, the fibre sequence is set-up for homeomorphisms, and under restrictions on the open manifolds under consideration. While we believe one could compare the two fibre sequences (which for instance requires setting up Cardim's fibre sequence in the smooth category, or our fibre sequence in the topological category), we do not carry out any such comparison.
\subsection{Pro-ends\protect\footnote{The contents of this section are not due to the sole work of the author. Jan Steinebrunner suggested the main strategy of this section, to whom we owe most credit. While carrying out the suggested strategy, the author discussed it with Thomas Blom, who suggested that Theorem \ref{thmE} ought to hold; finally, Maxime Ramzi provided a proof thereof.}}
We now turn to the second part of Objective (i), namely writing down a functor that sends a manifold to its pro-end. We begin with an informal discussion. Let $M$ be a manifold, and fix a choice of a compact exhaustion $\{K_i\}$ of $M$; denote the consecutive bordisms by $W_{i,i+1}\colon\partial K_i\rightsquigarrow\partial K_{i+1}$. Then, it follows immediately that for each $i$, we have an equivalence
$$M\setminus K_i \simeq W_{i,i+1}\cup_{\partial K_{i+1}} W_{i+1,i+2}\cup_{\partial K_{i+2}}W_{i+2,i+3}\cup_{\partial K_{i+3}}\cdots$$
where the latter denotes the iterated pushouts. Recall that pushouts are the compositions in $\Cosp(\catsingle{S})$; this suggests that one may define a functor
$$\tau\colon\PSh^\mathbb{N}(\Cosp(\catsingle{S})) \to \Cosp(\Pro(\catsingle{S}))$$
which, on objects, glues all the bordisms above degree $i$ to obtain the complement of $K_i$. In order to write down the functor $\tau$, the following proved necessary.
\begin{bigthm}\label{thmE}
    Fix an $\infty$-category $\catsingle{C}$. The functor $\Phi\colon\Fun(\mathbb{N},\catsingle{C})\to \PSh^\mathbb{N}(\catsingle{C})$, sending a sequential diagram in $\catsingle{C}$ to its colimit in $\PSh(\catsingle{C})$, descends to an equivalence $$\Fun(\mathbb{N},\catsingle{C})[W^{-1}] \xrightarrow[]{\simeq}\PSh^\mathbb{N}(\catsingle{C})$$
    where $W$ is the collection of cofinal re-indexing morphisms in $\Fun(\mathbb{N},\catsingle{C})$.
\end{bigthm}
With Theorem \ref{thmE} at hand, we proceed as follows: one may define a "tail" functor
$$\mathrm{tail}\colon \Fun(\mathbb{N},\Cosp(\catsingle{S})) \to \Cosp(\Pro(\catsingle{S}))$$
implementing the intuition above, namely it glues all the cospans above degree $i$ to obtain the complement. It can furtermore be shown that this functor sends all cofinal re-indexings to equivalences; see §\ref{section:tail-functor}. By Theorem \ref{thmE}, this descends to a functor on $\PSh^\mathbb{N}(\Cosp(\catsingle{S}))$, which we denote
$$\tau\colon\PSh^\mathbb{N}(\Cosp(\catsingle{S}))\to \Cosp(\Pro(\catsingle{S}))$$
We furthermore have a functor $u\colon\BordInf_d \to \Cosp(\catsingle{S})$, which on morphisms sends a bordism to the underlying cospan in spaces. We thus define the functor $\proEnd$ as the following composite:
\begin{center}
    \begin{tikzcd}
        \proEnd\colon\ManInf_d \arrow{r}{\bordEnd} & \PSh^\mathbb{N}(\BordInf_d)\arrow{r}{\PSh^\mathbb{N}(u)} & \PSh^\mathbb{N}(\Cosp(\catsingle{S})) \arrow{r}{\tau} & \Cosp(\Pro(\catsingle{S})) 
    \end{tikzcd}
\end{center}
By construction, this functor satisfies the required property, namely that if $M\in \OpenMan_d$ and $\{K_i\}$ is a chosen compact exhaustion of $M$, $\proEnd(M)\simeq\{M\setminus K_i\}$ as objects in $\Pro(\catsingle{S})$; furthermore, by construction, the triangle in Theorem \ref{thmA} commutes, where the dashed arrow in the theorem is the composite $\tau\circ \PSh^\mathbb{N}(u)$.

\subsection*{Acknowledgments}
I would like to thank Søren Galatius and Manuel Krannich for the many helpful discussions that led to this work, and Branko Juran for general discussions about free colimit completions; I would also like to thank Branko Juran and Manuel Krannich for reading a previous draft of this manuscript. My gratitude extends to Thomas Blom, Maxime Ramzi and Jan Steinebrunner for their help regarding the contents of §\ref{section:pro-ends}: the general strategy is due to Jan Steinebrunner. Thomas Blom pointed out that Theorem \ref{thm:localisation} ought to be true, which was a crucial input for the construction of the functor in §\ref{section:pro-ends}, and Maxime Ramzi provided a proof thereof.

During this work, I was partially supported by the Danish National Research Foundation through the Copenhagen Center for Geometry and Topology (DNRF151), and the European Union through an ERC grant (MaFC, 101221003).

\subsection*{AI disclosure}
During the preparation of this paper, the author used Claude Opus 5.5 (Anthropic), for the following tasks: checking the homotopy-limit computation underlying Lemma \ref{lem:invertible-emb-diff}(ii);  locating literature on topologies on groups of compactly supported diffeomorphisms; and proofreading the manuscript. In every case, the questions put to the AI concerned material the author had already worked out independently. The AI served as a checking tool; no result, idea, proof or construction in this paper originates from it. The entirety of the paper was written by the author, and the author takes full responsibility for all the content.

\section{Categorical preliminaries}
This section is a collection of $\infty$-categorical facts that are used throughout in this work.

\subsection{Definition and universal property}\label{sec:definition}
 The main references for free colimit completions are \cite[§5.3.6]{HTT}, \cite{Rezk}. Let $\mathcal{F}\subset \Cat_\infty$ be a class of small $\infty$-categories. The free $\mathcal{F}$-colimit completion of a category $\catsingle{C}$ is an $\infty$-category $\PSh^\mathcal{F}(\catsingle{C})$, uniquely characterized by the following:
\begin{enumerate}
    \item $\PSh^\mathcal{F}(\catsingle{C})$ is a full subcategory of $\PSh(\catsingle{C})$ containing the image of the Yoneda embedding $\catsingle{C} \hookrightarrow \PSh(\catsingle{C})$;
    \item The category $\PSh^\mathcal{F}(\catsingle{C})$ is stable under $\mathcal{F}$-colimits in $\PSh(\catsingle{C})$;
    \item This category satisfies the following universal property: for any $\infty$-category $\catsingle{D}$ which admits $\mathcal{F}$-indexed colimits, precomposition along the Yoneda embedding yields an equivalence of categories
    $$\Fun^\mathcal{F}(\PSh^\mathcal{F}(\catsingle{C}),\catsingle{D}) \to \Fun(\catsingle{C},\catsingle{D})$$
    where $\Fun^\mathcal{F}$ denotes the category of $\mathcal{F}$-colimit preserving functors. In more concrete terms, given a category $\catsingle{D}$ which admits $\mathcal{F}$-colimits, any functor $\catsingle{C} \to \catsingle{D}$ admits an essentially unique extension to a functor $\PSh^\mathcal{F}(\catsingle{C}) \to \catsingle{D}$.
\end{enumerate}
One may replace (iii) in the above with the following condition, which will be very useful in practice:
\begin{enumerate}
    \item[(iii)'] $\PSh^\mathcal{F}(\catsingle{C})$ is the smallest full subcategory of $\PSh(\catsingle{C})$ satisfying (i) and (ii).
\end{enumerate}

\subsubsection{Free sequential colimit completions}

A particular choice of a class of small $\infty$-categories will play a crucial role in this work, namely $\mathcal{F}=\{\mathbb{N}\}$, where $\mathbb{N}$ refers to the poset of natural integers, viewed as a category. We restrict to this case here, and investigate some of the special features this category has, which may not hold for arbitrary $\mathcal{F}$.
\begin{dfn}
    Fix $\mathcal{F}$ a collection of small $\infty$-categories and an $\infty$-category $\catsingle{C}$. We let $\PSh^\mathcal{F}_{nv}(\catsingle{C})$ be the full subcategory of $\PSh(\catsingle{C})$ on those objects which are $\mathcal{F}$-colimits of representables in $\PSh(\catsingle{C})$.
\end{dfn}
It is clear from the definitions that one obtains an inclusion $\PSh^\mathcal{F}_{nv}(\catsingle{C}) \subseteq \PSh^\mathcal{F}(\catsingle{C})$; this is because $\PSh^\mathcal{F}(\catsingle{C})$ contains all representable presheaves, and is closed under $\mathcal{F}$-colimits in $\PSh(\catsingle{C})$. However, it is not clear in general that these two categories are equivalent. Indeed, it is not clear whether $\PSh^\mathcal{F}_{nv}(\catsingle{C})$ is itself stable under $\mathcal{F}$-colimits in $\PSh(\catsingle{C})$. This, however, does hold when $\mathcal{F}=\{\mathbb{N}\}$, as we will show in a moment. However, we begin by observing the following description of the mapping space of $\PSh^\mathbb{N}_{nv}(\catsingle{C})$, which will be very useful for us. This follows immediately from \cite[Proposition 2.1]{efimov}.
\begin{lem}\label{lem:interchange-lemma}
    Let $\mathcal{G}$ denote the directed poset of strictly monotone maps $f\colon \mathbb{N} \to \mathbb{N}$, ordered by $f\geq g$ if $f(n) \geq g(n)$ for all $n\in \mathbb{N}$. Let $\{c_i\}$, $\{d_j\}$ be objects in $\Fun(\mathbb{N},\catsingle{C})$, where $\catsingle{C}$ is an $\infty$-category, and denote by $\colim_i h_{c_i}$ and $\colim_j h_{d_j}$ their respective colimits in $\PSh(\catsingle{C})$, where $h_c$ denotes the representable presheaf on $c$. Then, there is a natural equivalence
    $$\Map_{\PSh(\catsingle{C})}(\colim_{i \in \mathbb{N}} h_{c_i}, \colim_{j\in \mathbb{N}} h_{d_j}) \simeq \colim_{f\in \mathcal{G}} \Map_{\Fun(\mathbb{N},\catsingle{C})}(\{c_i\},f^\ast\{d_{i}\})$$
\end{lem}
\begin{proof}
    We begin with the equivalence
    $$\Map_{\PSh(\catsingle{C})}(\colim_i h_{c_i},\colim_jh_{d_j}) \simeq \lim_i \colim_j \Map_{\catsingle{C}}(c_i,d_j)$$
    By \cite[Proposition 2.1]{efimov}, one may interchange the $\lim$-$\colim$ in the above to obtain a natural equivalence
    $$\lim_i \colim_j \Map_{\catsingle{C}}(c_i,d_j) \simeq \colim_{\varphi\colon\mathbb{N} \to \mathbb{N}} \lim_{i\leq j} \Map_{\catsingle{C}}(c_i,d_{\varphi(j)})$$
    where the limit is taken over the twisted arrow category of $\mathbb{N}$, and where $\varphi$ ranges over monotone self-maps of $\mathbb{N}$. For a fixed $\varphi$, we have an equivalence
    $$\lim_{(i\leq j)\in \tw(\mathbb{N})}\Map_{\catsingle{C}}(c_i,d_{\varphi(j)}) \simeq \Map_{\Fun(\mathbb{N},\catsingle{C})}(\{c_i\},\varphi^\ast \{d_j\})$$
    as follows from \cite[Proposition 5.1]{GepnerHaugsengNikolaus}. As the inclusion of $\mathcal{G}$ into the set of monotone maps $\mathbb{N} \to \mathbb{N}$ is clearly cofinal, the result follows.
\end{proof}
We now show that $\PSh^\mathbb{N}_{nv}(\catsingle{C})$ is equivalent to $\PSh^\mathbb{N}(\catsingle{C})$.
\begin{lem}\label{lem:naive-colim}
    Let $\catsingle{C}$ be an $\infty$-category. The full subcategory inclusion $\PSh^\mathbb{N}_{nv}(\catsingle{C}) \hookrightarrow \PSh(\catsingle{C})$ is closed under sequential colimits. As a consequence, $\PSh^\mathbb{N}_{nv}(\catsingle{C})$ is equivalent to $\PSh^\mathbb{N}(\catsingle{C})$.
\end{lem}
\begin{proof}
    Consider a diagram $X_\bullet \colon \mathbb{N} \to \PSh^\mathbb{N}_{nv}(\catsingle{C})$. For each $n \in \mathbb{N}$, since $X_n \in \PSh^\mathbb{N}_{nv}(\catsingle{C})$, we note that
    $$X_n \simeq \colim_m Y_{n,m}$$
    where $Y_{n,m}$ is a representable presheaf. Now, as $\PSh(\catsingle{C})$ admits colimits, we let
    $$X \coloneqq \colim_n X_n$$
    Now, we investigate the maps $\varphi_{i,i+1}\colon X_i \to X_{i+1}$. By Lemma \ref{lem:interchange-lemma}, $\varphi_{i,i+1}$ is equivalent to the following: a choice of a monotone map $f_{i,i+1} \colon \mathbb{N} \to \mathbb{N}$, a natural transformation $\alpha_{i,i+1} \colon \{Y_{i,n}\} \to \{Y_{i+1,f(n)}\}$ in $\Fun(\mathbb{N},\catsingle{C})$, and an equivalence $\varphi_{i,i+1}\simeq \colim(\alpha_i)$, where $\colim(\alpha_f)$ is the induced map on the colimits. We thus obtain a sequence of natural transformations
    $$\{Y_{0,n}\} \xrightarrow[]{\alpha_{0,1}} \{Y_{1,f_{0,1}(n)}\} \xrightarrow{\alpha_{1,2}}Y_{2,f_{1,2}\circ f_{0,1}(n)} \xrightarrow[]{\alpha_{2,3}} \cdots$$
    in $\Fun(\mathbb{N},\catsingle{C})$, i.e. a functor $\mathbb{N} \to \Fun(\mathbb{N},\catsingle{C})$, which upon applying colimits, is equivalent to the original sequence
    $$X_0 \xrightarrow[]{\varphi_{0,1}}X_1 \xrightarrow[]{\varphi_{1,2}}X_2 \xrightarrow[]{\varphi_{2,3}}\cdots$$
    as elements in $\Fun(\mathbb{N},\PSh^\mathbb{N}(\catsingle{C}))$; here, we are using that $\mathbb{N}$ is the iterated pushouts of the walking arrow category, and hence prescribing an equivalence between two such sequences only requires us to provide homotopies between the squares (and no higher homotopies). By a diagonal argument, the statement follows.
\end{proof}

\begin{rem}
    We comment on the functoriality of $\PSh^\mathbb{N}(-)$. Given two $\infty$-categories $\catsingle{C}$, $\catsingle{D}$ and a functor $f \colon \catsingle{C} \to \catsingle{D}$, we claim that the functor $\PSh(f)\colon\PSh(\catsingle{C})\to \PSh(\catsingle{D})$ restricts to a functor on $\PSh^\mathbb{N}$. More precisely, the functoriality of $\PSh(-)$ yields the leftmost commutative diagram in the below diagram
    \begin{center}
        \begin{tikzcd}
            \catsingle{C} \arrow[hook]{r} \arrow{d}[swap]{f} & \PSh(\catsingle{C}) \arrow{d}{\PSh(f)} & \arrow[hook]{l}\PSh^\mathbb{N}(\catsingle{C}) \arrow[dashed]{d}
            \\
            \catsingle{D} \arrow[hook]{r} & \PSh(\catsingle{D}) & \PSh^\mathbb{N}(\catsingle{D}) \arrow[hook]{l}
        \end{tikzcd}
    \end{center}
    As indicated in the diagram, we claim that $\PSh(f)$, when restricted to $\PSh^\mathbb{N}(\catsingle{C})$, lands in $\PSh^\mathbb{N}(\catsingle{D})$. As $\PSh^\mathbb{N}(\catsingle{D})$ is a full subcategory of $\PSh(\catsingle{D})$, this can be checked on objects; to that end, we fix $X \in \PSh^\mathbb{N}(\catsingle{C})$. By Lemma \ref{lem:naive-colim}, it follows that $X\simeq \colim_{i\in\mathbb{N}}h_{c_i}$, where $h_{c_i}$ denotes the representable presheaf on objects $c_i \in \catsingle{C}$. We have the following string of equivalences
    \begin{align*}
        \PSh(f)(X) \simeq \PSh(f)(\colim_{i\in\mathbb{N}}h_{c_i})\simeq \colim_{i\in\mathbb{N}}\PSh(f)(h_{c_i})\simeq \colim_{i\in\mathbb{N}}h_{f(c_i)}
    \end{align*}
    where again, $h_{f(c_i)}$ denotes the representable presheaf in $\PSh(\catsingle{D})$ on the objects $f(c_i)$. The above equivalences follow from the obvious facts that $\PSh(f)$ preserves both colimits and representable objects. We thus denote the restriction of $\PSh(f)$ to $\PSh^\mathbb{N}$ by $\PSh^\mathbb{N}(f)$.
\end{rem}
We finally list a lemma that will be useful in the following. This gives a list of necessary and sufficient conditions of a functor $\PSh^\mathbb{N}(\catsingle{C}) \to \catsingle{D}$ to be an equivalence. Its proof is a straightforward adaption of the proof of \cite[5.3.5.11]{HTT}, which we therefore omit.
\begin{lem}\label{lem:conditions}
    Let $\catsingle{C}$ be an $\infty$-category, and let $\catsingle{D}$ be an $\infty$-category admitting sequential colimits. Let $F \colon \PSh^\mathbb{N}(\catsingle{C}) \to \catsingle{D}$ be a sequential colimit-preserving functor. $F$ is an equivalence if and only if
    \begin{enumerate}
        \item $f\coloneqq F \circ y$, where $y$ is the Yoneda embedding, is fully faithful;
        \item For each object $c\in \catsingle{C}$, $f(c)$ is sequentially compact;
        \item The category $\catsingle{D}$ is generated under sequential colimits by the image of $f$.    \end{enumerate}
\end{lem}

As a direct application of the above, we show the following concerning the behaviour of slice categories under $\PSh^\mathbb{N}$.
\begin{lem}\label{lem:slice}
    Let $\catsingle{C}$ be an $\infty$-category, and fix $c \in \catsingle{C}$ an object. Let $y \colon \catsingle{C} \hookrightarrow \PSh^\mathbb{N}(\catsingle{C})$ denote the Yoneda embedding.
    \begin{enumerate}
        \item There is an equivalence of categories
        $$\PSh^\mathbb{N}(\catsingle{C}_{c/}) \xrightarrow[]{\simeq} \PSh^\mathbb{N}(\catsingle{C})_{y(c)/}$$
        \item There is an equivalence of categories
        $$\PSh^\mathbb{N}(\catsingle{C}_{/c}) \xrightarrow[]{\simeq} \PSh^\mathbb{N}(\catsingle{C})_{/y(c)}$$
    \end{enumerate}
\end{lem}
\begin{proof}
    We begin by showing (i). The Yoneda embedding $y \colon \catsingle{C} \hookrightarrow \PSh^\mathbb{N}(\catsingle{C})$ induces a functor
    $$y_{c/}\colon\catsingle{C}_{c/} \to \PSh^\mathbb{N}(\catsingle{C})_{y(c)/}$$
    We observe that $\PSh^\mathbb{N}(\catsingle{C})_{y(c)/}$ admits sequential colimits. This follows from combining the two following facts:
    \begin{itemize}
        \item $\PSh^\mathbb{N}(\catsingle{C})$ admits sequential colimits by definition;
        \item For any category $\catsingle{D}$ and $d \in \catsingle{D}$, the forgetful functor $\catsingle{D}_{d/} \to \catsingle{D}$ creates all contractible colimits that $\catsingle{D}$ admits \cite[4.4.2.9]{HTT} (where the forgetful functor $\catsingle{D}_{d/} \to \catsingle{D}$ is a left fibration by \cite[2.1.2.2]{HTT}).
    \end{itemize}
    By the universal property of $\PSh^\mathbb{N}$, we obtain a unique sequential colimit preserving functor $\varphi$ making the following diagram commute:
    \begin{center}
        \begin{tikzcd}
            \catsingle{C}_{c/} \arrow{r}{y_{c/}} \arrow[hook]{d} & \PSh^\mathbb{N}(\catsingle{C})_{y(c)/}
            \\
            \PSh^\mathbb{N}(\catsingle{C}_{c/}) \arrow[dashed]{ur}[swap]{\varphi}
        \end{tikzcd}
    \end{center}
    We now make use of Lemma \ref{lem:conditions} to show that $\varphi$ is an equivalence.
    \begin{enumerate}
        \item We first observe that $y_{c/}$ is a fully faithful functor. Indeed, the Yoneda embedding $y \colon \catsingle{C} \to \PSh^\mathbb{N}(\catsingle{C})$ is fully faithful, and $y_{c/}$ is obtained from $y$ by slicing in the source; this clearly remains fully faithful.
        \item We need to show that for every $(c \xrightarrow[]{f}x) \in \catsingle{C}_{c/}$, $y_{c/}(c \xrightarrow[]{f}x)$ is a sequentially compact object in $\PSh^\mathbb{N}(\catsingle{C})_{y(c)/}$. Let $\{(c \xrightarrow[]{g_i})y_i\}_{i \in \mathbb{N}}$ be a sequential diagram in $\PSh^\mathbb{N}(\catsingle{C})_{y(c)/}$; we observe that the colimit of the above diagram is mapped, by the forgetful functor, to the colimits of $y_i$ in $\PSh^\mathbb{N}(\catsingle{C})$. Consider the following commutative diagram, where for the purpose of clarity, we replace the notation $c \to x$ simply by $x$
        \begin{center}
            \begin{tikzcd}
                \Map_{\PSh^\mathbb{N}(\catsingle{C})_{c/}}(x, \colim_i y_i) \arrow{r}  & \Map_{\PSh^\mathbb{N}(\catsingle{C})}(x,\colim_i y_i) \arrow{r}{f^\ast} \arrow{d}{\simeq}[swap]{\circled{1}} & \Map_{\PSh^\mathbb{N}(\catsingle{C})}(c,\colim_i y_i) \arrow{d}{\simeq}[swap]{\circled{2}}
                \\
                \fib_g(\colim_i (f^\ast)) \arrow{r} \arrow[no head]{d}{\simeq}[swap]{\circled{3}} & \colim_i \Map_{\PSh^\mathbb{N}(\catsingle{C})}(x,y_i) \arrow{r}[swap]{\colim_i f^\ast} & \colim_i \Map_{\PSh^\mathbb{N}(\catsingle{C})}(c,y_i)
                \\
                \colim_i \Map_{\PSh^\mathbb{N}(\catsingle{C})_{y(c)/}}(x,Y_i)
            \end{tikzcd}
        \end{center}
        where
        \begin{itemize}
            \item Each row in the above square is a fibre sequence (with the top row exhibiting the mapping space in the slice category as the fibre at $g$);
            \item $\circled{1}$ and $\circled{2}$ are equivalences as any representable in $\PSh^\mathbb{N}(\catsingle{C})$ is sequentially compact;
            \item \circled{3} is an equivalence, as finite limits commute with filtered colimits, of which sequential colimits is an example.
        \end{itemize}
        \item This step requires us to show that any object in $\PSh^\mathbb{N}(\catsingle{C})_{c/}$ is generated under sequential colimits from objects in $\catsingle{C}_{c/} \subset \PSh^\mathbb{N}(\catsingle{C}_{c/})$. Fix such an object $c \to X$ where $X \in \PSh^\mathbb{N}(\catsingle{C})$. Applying the forgetful functor, we observe that $X \simeq \colim_i x_i$, with $x_i \in \catsingle{C}$, by Lemma \ref{lem:naive-colim}. As the forgetful functor creates contractible colimits, it follows that this colimit diagram lifts uniquely to a colimit diagram in $\PSh^\mathbb{N}(\catsingle{C})_{c/}$, and in particular $X \simeq \colim_i (c\to x_i)$. 
    \end{enumerate}
    As a consequence of Lemma \ref{lem:conditions}, it follows that $\varphi$ is an equivalence.
    \\

    To show (ii), this follows from the same strategy; in fact, it is more simple, as colimits in $\catsingle{C}_{/c}$ are computed underlying.
\end{proof}
We finally comment on the right adjoint of the free sequential colimit completion. 
\begin{lem}\label{lem:right-adjoint-colim-extension}
    Let $\catsingle{C}$ and $\catsingle{D}$ be two $\infty$-categories, and assume that $\catsingle{D}$ admits sequential colimits. Let $f \colon \catsingle{C} \to \catsingle{D}$ be a functor, $f_!$ its unique sequential-colimit preserving extension, and assume that $f$ has the property that $f^\ast \colon \PSh(\catsingle{D}) \to \PSh(\catsingle{C})$ restricts to a functor $f^\ast \colon \PSh^\mathbb{N}(\catsingle{D}) \to \PSh^\mathbb{N}(\catsingle{C})$. Then, $f_!$ is left adjoint to $f^\ast h \colon \catsingle{D} \to \PSh^\mathbb{N}(\catsingle{C})$, where $h$ denotes the Yoneda embedding. 
\end{lem}
\begin{proof}
    We need to show that, for all $d \in \catsingle{D}$, we have a natural equivalence
    $$\Map_{\catsingle{D}}(f_!(-),d) \simeq \Map_{\PSh^\mathbb{N}(\catsingle{C})}(-,f^\ast h(d))$$
    as functors $\PSh^\mathbb{N}(\catsingle{C})^\op \to \catsingle{S}$. We first observe that $f_!(-)$ preserves sequential colimits, and therefore both the above functors send sequential colimits to limits over the same diagram. Thus, it suffices to check that the precompositions of the functors with the Yoneda embedding yield equivalent functors. We consider the following equivalences
    \begin{align*}
        \Map_{\PSh^\mathbb{N}(\catsingle{C})}(h(-),f^\ast h(d)) \simeq \Map_{\PSh(\catsingle{C})}(h(-),f^\ast h(d))&\simeq \Map_{\PSh(\catsingle{D})}(f_! h(-),h(d))
        \\
        &\simeq \Map_{\catsingle{D}}(f(-),d)
    \end{align*}
    and all the above equivalences are natural.
\end{proof}
 \subsection{$\PSh^\mathbb{N}$ as a Dwyer-Kan localization of $\Fun(\mathbb{N},-)$}
    For $\catsingle{C}$ an $\infty$-category, taking colimits in $\PSh(\catsingle{C})$ defines a functor
    $$\Phi \colon \Fun(\mathbb{N},\catsingle{C}) \to \PSh(\catsingle{C})$$
    which lands in $\PSh^\mathbb{N}(\catsingle{C})$. We claim that the functor $\Phi$ above exhibits $\PSh^\mathbb{N}(\catsingle{C})$ as the Dwyer-Kan localization of $\Fun(\mathbb{N},\catsingle{C})$ at a class of morphisms we describe now.
    \\

    Consider the poset of monotone self-maps $f \colon \mathbb{N} \to \mathbb{N}$, where the partial order is given by $g \leq g'$ if and only if $g(n) \leq g'(n)$ for all $n \in \mathbb{N}$; let $\mathcal{G}$ be the subposet on those strictly monotone maps $f \colon \mathbb{N} \to \mathbb{N}$. For $F \in \Fun(\mathbb{N}, \catsingle{C})$, precomposition with an element $g \in \mathcal{G}$ defines a natural transformation
    $$\eta_g^F \colon F \to F \circ g$$
    We set $W \coloneqq \{\eta_g^F\,: F\in \Fun(\mathbb{N},\catsingle{C}),g \in \mathcal{G}\}$, and call them the \emph{cofinal re-indexing morphisms in $\Fun(\mathbb{N},\catsingle{C})$}. It is clear that the functor $\colim$ inverts all reindexing maps, which in particular induces a map
    $$\widetilde{\Phi}\colon\Fun(\mathbb{N},\catsingle{C})[W^{-1}] \to \PSh^\mathbb{N}(\catsingle{C})$$
    where the source is the Dwyer-Kan localization of $\Fun(\mathbb{N},\catsingle{C})$ at $W$. This functor is an equivalence, and exhibits $\PSh^\mathbb{N}(\catsingle{C})$ as the Dwyer-Kan localization of $\Fun(\mathbb{N},\catsingle{C})$ at the class of the cofinal re-indexing morphisms; we follow the conventions and definitions of \cite{derivedmappingspacesinftycategories} for Dwyer-Kan localizations.
    \begin{thm}\label{thm:localisation}
        The functor $\widetilde{\Phi}\colon \Fun(\mathbb{N},\catsingle{C})[W^{-1}] \to \PSh^\mathbb{N}(\catsingle{C})$ is an equivalence of categories.
    \end{thm}
    \begin{proof}
    That this functor is essentially surjective follows readily from Lemma \ref{lem:naive-colim}. The heart of the argument is to show that this functor is fully faithful.

    We set up some notation. Let $\mathcal{L}\coloneqq \Fun(\mathbb{N},\catsingle{C})[W^{-1}]$, and let $\gamma\colon\Fun(\mathbb{N},\catsingle{C}) \to \mathcal{L}$ be the localization functor. By \cite[Proposition 1.3]{derivedmappingspacesinftycategories}, the restriction functor
    $$\gamma^\ast\colon \PSh(\mathcal{L}) \to \PSh(\Fun(\mathbb{N},\catsingle{C}))$$
    is fully faithful, with image the $W$-local presheaves, i.e. functors $(\Fun(\mathbb{N},\catsingle{C}))^\op \to \catsingle{S}$ sending the morphisms in $W$ to equivalences in $\catsingle{S}$. For $d_\bullet \in \Fun(\mathbb{N},\catsingle{C})$, let $h_d$ denote the representable functor in $\PSh(\Fun(\mathbb{N},\catsingle{C}))$, and let $G_d$ denote the functor $\Map_{\PSh^\mathbb{N}(\catsingle{C})}(\Phi(-),\Phi(d))$. We begin by observing that $G_d$ is $W$-local: indeed, if $f\in \mathcal{G}$, and for any $c_\bullet \in \Fun(\mathbb{N},\catsingle{C})$ then $\Phi(c) \to \Phi(f^\ast c)$ is by definition the comparison map
    $$\colim_{i\in\mathbb{N}}c_i \to \colim_{i\in\mathbb{N}}c_{f(i)}$$
    which is an equivalence, as $f$ is a cofinal re-indexing. Furthermore, the functor $\Phi$ provides a map $\Phi\colon h_d \to G_d$, which factors uniquely into a diagram
    \begin{center}
        \begin{tikzcd}
            h_d \arrow{r}{\Phi} \arrow{d} & G_d
            \\
            \gamma^\ast \gamma_! h_d \arrow[dashed]{ur}{\exists!}
        \end{tikzcd}
    \end{center}
    where the vertical map is induced by the unit of the $(\gamma_!,\gamma^\ast)$-adjunction; we denote the dashed map by $\widetilde{\Phi}$. Note that for any $c \in \mathcal{L}$, where we may assume $c$ is of the form $\gamma(c)$ for some $c \in \Fun(\mathbb{N},\catsingle{C})$, the map $\widetilde{\Phi}(\gamma(c))$ yields the map
    $$\Map_\mathcal{L}(\gamma(c),\gamma(d)) \to \Map_{\PSh^\mathbb{N}(\catsingle{C})}(\widetilde{\Phi}\gamma(c),\widetilde{\Phi}\gamma(d))$$
    which is indeed the map induced by the functor $\widetilde{\Phi}\colon\mathcal{L} \to \PSh^\mathbb{N}(\catsingle{C})$ on mapping spaces. Consequently, fully-faithfulness of $\widetilde{\Phi}$ is now equivalent to the map $\widetilde{\Phi}\colon \gamma^\ast \gamma_! h_d \to G_d$ being an equivalence in $\PSh(\mathcal{L})$, for all $d \in \Fun(\mathbb{N},\catsingle{C})$.

    We now exhibit $G_d$ as a particular colimit of representables. By Lemma \ref{lem:interchange-lemma}, we obtain an equivalence
    $$G_d \simeq \colim_{f\in\mathcal{G}}h_{f^\ast d}$$
    In particular, for all $d\in \Fun(\mathbb{N},\catsingle{C})$, the functor $\mathcal{G} \to \Fun(\mathbb{N},\catsingle{C})_{d/}$, sending $f$ to the pair $(d\to f^\ast d)$ is a $\mathcal{G}$-resolution of $d$, following \cite[Definition 2.1]{derivedmappingspacesinftycategories} (the other two conditions in this definition are easily verified in our case), and therefore, by Theorem 2.2 in \emph{loc.cit.}, it follows that
    $$\gamma^\ast \gamma_! h_d \to G_d$$
    is an equivalence, which implies the result.
    \end{proof}
\bigskip

\section{The category of manifolds as a formal completion of compact manifolds}\label{section:OpenManifolds}
We first give an informal description of the categories of manifolds that are of interest in this work; we shall return in a moment to give a more detailed definition. Let $\ManInf_d$ denote the $\infty$-category of smooth $d$-dimensional manifolds with empty boundary, and spaces of smooth embeddings. Let $\CpctMan_d$ denote the category of smooth, compact $d$-dimensional manifolds, and spaces of smooth embeddings between them. There is a functor $\interior\colon \CpctMan_d \to \ManInf_d$, sending a compact $d$-manifold to its interior, and restricting an embedding between two manifolds to the interior; as any manifold is isotopy equivalent to its interior, the above functor is fully faithful. Additionally, we let $\BordInf_d$ be the $\infty$-category of bordisms: that is, this category has as space of objects all closed $(d-1)$-dimensional smooth manifolds, and as space of morphisms the space of compact smooth manifolds with boundary. We shall return in a moment to give a rigorous definition of all these categories; but first, we already state the main result we work towards in this section. 

\begin{thm}\label{thm:openMfldFreeColim} 
    The unique sequential colimit preserving extension of the functor $\interior\colon \CpctMan_d  \hookrightarrow \OpenMan_d $ induces an equivalence of $\infty$-categories
    $$\interior_! \colon \PSh^\mathbb{N}(\CpctMan _d)  \xrightarrow[]{\simeq} \OpenMan_d $$
    Its inverse admits a description as a restricted Yoneda:
    $$\Exh\colon \OpenMan_d  \hookrightarrow \PSh^\mathbb{N}(\OpenMan_d  ) \xrightarrow[]{\interior^\ast} \PSh^\mathbb{N}(\CpctMan_d )$$   
\end{thm}
It is to be pointed out that there are two implicit claims in the above theorem statement: first, that there exists a unique sequential colimit preserving extension, which requires the category $\OpenMan_d $ to admit sequential colimits; second, that the pullback functor $\interior^\ast$ exists. Both these statements are addressed in the following sections, as Lemmas \ref{lem:OpenMan-sequential-colimits} and \ref{lem:restriction-functor-exists} respectively.

One upshot of the above theorem is that one can construct a functor $\PSh^{\mathbb{N}}(-)\colon \CatInf \to \CatInf$ defined via a universal property, which sends $\CpctMan_d$ to $\OpenMan_d$; that is, one recovers the category of all manifolds formally from that of compact manifolds.

In the following subsections, we define all the objects invoked in the above theorem, and conclude with a proof.
\subsection{Categories of manifolds}
We begin by describing the categories $\ManInf_d$ and $\CpctMan_d$ in more details. Let $\mathrm{Mfld}_d$ denote the topologically enriched category, where objects are smooth $d$-dimensional manifolds with empty boundary (viewed as a subset of $\mathbb{R}^\infty$), and where the morphisms are the topological spaces of smooth embeddings, topologized via the weak Whitney $C^\infty$-topology (as in \cite[Chapter 2]{Hirsch}). Let $\mathrm{Mfld}^{\mathrm{cpt}}_d$ denote the topologically enriched category of smooth compact $d$-dimensional manifolds (potentially with boundary), and spaces of embeddings between them, topologized via the (weak) Whitney $C^\infty$-topology. We denote $\ManInf_d$ (resp. $\CpctMan_d$) as the $\infty$-category obtained from $\mathrm{Mfld}_d$ (resp. $\mathrm{Mfld}^{\mathrm{cpt}}_d$) by applying coherent nerves. There is a functor $\interior \colon \mathrm{Mfld}^{\mathrm{cpt}}_d \to \mathrm{Mfld}_d$, obtained by sending a compact manifold to its interior, and restricting an embedding between two compact manifolds to the embedding of their interiors. By abuse of notation, we let $\interior\colon \CpctMan_d \to \ManInf_d$ denote the induced functor between the $\infty$-categories. Our first point of order is to show that this functor is fully faithful.
\begin{lem}\label{lem:interior-functor-fully-faithful}
    Let $M$ and $N$ be two smooth $d$-manifolds. Then, the map
    $$\Emb(M,N) \to \Emb(\interior(M), \interior(N))$$
    given by the restriction of an embedding $f$ to the interior of $M$, is a homotopy equivalence.
\end{lem}
\begin{proof}
    Choosing a boundary collar for both $M$ and $N$, we observe that the embeddings $\interior(M) \hookrightarrow M$ and $\interior(N) \hookrightarrow N$ both admit isotopy inverses $M \hookrightarrow \interior(M)$ and $N \hookrightarrow \interior(N)$. The result then follows readily.
\end{proof}

In particular, one observes that after coherent nerve, the functor $\interior\colon\CpctMan_d \to \OpenMan_d$ is fully faithful.
\begin{rem}
    The fully faithfulness of $\interior$ identifies, in particular, the category $\CpctMan_d$, with the full subcategory of $\ManInf_d$ of those manifolds that happen to be isotopy equivalent (hence diffeomorphic, see Lemma \ref{lem:invertible-emb-diff} below) to the interior of a compact manifold. We will be using interchangeably this identification. Indeed, at some later points, one may see a statement along the lines of "let $K$ be a compact manifold" next to a statement "$K \in \ManInf_d$"; by this, we are identifying $K$ with $\interior(K)$, in order to declutter notation.
\end{rem}

Given a smooth $d$-manifold $W$, a choice of a proper Morse function $f \colon W \to \mathbb{R}$ \cite[Corollary 6.7]{MilnorMorse} yields an \emph{exhaustion of $W$ by compact submanifolds}: by this, we mean a sequence $\{K_i\}_{i\in\mathbb{N}}$ of compact, codimension 0 smooth submanifolds, such that $K_i \subseteq\interior(K_{i+1})$, and $W = \bigcup_i K_i$. As a sublation, what we now show is that the category $\OpenMan_d$ admits sequential colimits, and a choice of a compact exhaustion of a manifold $W$ exhibits $W$ as a sequential colimit diagram in $\OpenMan_d$.

\begin{lem}\label{lem:OpenMan-sequential-colimits}
    The category $\OpenMan_d$ admits sequential colimits.
\end{lem}
\begin{proof}
    Consider a functor $\mathbb{N} \to \OpenMan_d $. Viewing $\mathbb{N}$ as iterated pushouts of $[1]$ along $[0]$, this in particular gives rise to a sequence of arrows $W_i \to W_{i+1}$, together with isomorphisms between sources and targets. As morphisms in $\OpenMan_d $ are represented by embeddings, this can be strictified to a sequence
    $$W_1 \hookrightarrow W_2 \hookrightarrow W_3 \hookrightarrow \cdots$$
    of manifolds and embeddings, i.e. one does not need to keep track of coherences. Let $W \coloneqq \bigcup_i W_i$, endowed with the union topology. We observe that $W$ admits the structure of a smooth manifold. That it is a topological manifold is fairly straightforward, and we observe that smoothness can be checked on each individual $W_i$. We furthermore observe that we have a sequence of embeddings $W_i \hookrightarrow W$ in $\OpenMan_d $. It remains to show that for all $V \in \OpenMan_d $, the induced map
    $$\Emb (W,V) \to \lim_i \Emb (W_i,V)$$
    is an equivalence. This is an instance of the cofunctor $\Emb (-,V)$ being \emph{good} in the sense of Weiss \cite[Proposition 1.4]{Weiss}.
\end{proof}
\begin{rem}
    The following, while being straightforward, is a fact we will be using everywhere. Let $W$ be a manifold, and let $\{K_i\}_{i \in \mathbb{N}}$ be a compact exhaustion of $W$. Then, $W \simeq \colim_i K_i$ in $\OpenMan_d $, which follows rather immediately from the definition of the colimit in $\OpenMan_d $ above.
\end{rem}

The following crucial lemma identifies the notion of isotopy invertible embeddings and diffeomorphisms. We emphasize that item (i) below is due to Agol, as a Mathoverflow answer in \cite{Agol}.
\begin{lem}\label{lem:invertible-emb-diff}\,
\begin{enumerate}
    \item Fix $M$, $N$ two manifolds in $\OpenMan_d$, and let $f \colon M \hookrightarrow N$ be an embedding which admits an isotopy inverse. Then, $f$ is isotopic to a diffeomorphism.
    \item Let $W \in \OpenMan_d$. Let $\Emb^{\simeq}(W,W)$ denote the space of isotopy invertible self-embeddings of $W$. The map
    $$\Diff(W) \to \Emb^{\simeq}(W,W)$$
    is an equivalence of spaces.
\end{enumerate}
\end{lem}
\begin{proof}
    We first show (i), i.e. expand on the argument in \cite{Agol}. We fix an embedding $f\colon M \hookrightarrow N$ which is an isotopy equivalence, and we fix $g\colon N \hookrightarrow M$ its isotopy inverse. As these maps induce bijections on $\pi_0$, we restrict to the case where $M$ and $N$ are connected. Furthermore, as the statement is obvious for closed manifolds, we assume that both $M$ and $N$ are open. We consider the two following manifolds:
    \begin{align*}
        W_f &\coloneqq \colim\left(M \xrightarrow[]{f}N \xrightarrow[]{g}M \xrightarrow[]{f}N \xrightarrow[]{g}\cdots\right);
        \\
        W_g &\coloneqq \colim\left( N \xrightarrow[]{g}M \xrightarrow[]{f}N \xrightarrow[]{g}M \xrightarrow[]{f}\cdots \right) 
    \end{align*}
    where the colimit is computed in $\OpenMan_d$. A first point of observation is that we have a diffeomorphism $W_f \cong W_g$; indeed, we have a strictly commutative diagram
    \begin{center}
        \begin{tikzcd}
            M \arrow{r}{f} \arrow{d}[swap]{f} & N \arrow{r}{g} \arrow{d}[swap]{g} & M \arrow{r}{f} \arrow{d}[swap]{f} & N \arrow{r}{g} \arrow{d}[swap]{g} & \cdots
            \\
            N \arrow{r}[swap]{g} & M \arrow{r}[swap]{f} & N \arrow{r}[swap]{g} & M \arrow{r}[swap]{f} & \cdots
        \end{tikzcd}
    \end{center}
    which, upon taking the colimit in $\OpenMan_d$, yields an embedding $W_f \hookrightarrow W_g$; by the definition of the underlying topological space of the colimit in $\OpenMan_d$ (as the colimit in underlying topological spaces), it is clear that this map is a homeomorphism of underlying topological spaces, where one map is a smooth local diffeomorphism, hence in particular a diffeomorphism. We denote this diffeomorphism by $\varphi$.
    By definition of $\varphi$, we have a commutative square in $\OpenMan_d$, as follows:
    \begin{equation}\label{diagram:emb-v-diff}
        \begin{tikzcd}
            M \arrow[hook]{r}{f} \arrow[hook]{d} & M \arrow[hook]{d}
            \\
            W_f \arrow{r}{\cong}[swap]{\varphi} & W_g
        \end{tikzcd}
    \end{equation}
    where the vertical embeddings are obtained by embedding $M$ (resp. $N$) to the first factor of the colimit. We now argue that both $M\hookrightarrow W_f$ and $N \hookrightarrow W_g$ are isotopic to diffeomorphisms; this is an expansion of Agol's argument in \cite{Agol}. It will be sufficient to show it for the embedding $M \hookrightarrow W_f$, which we henceforth denote $\iota_0$. Let $H$ be a choice of an isotopy $g\circ f\simeq \id_M$, and choose a compact exhaustion of $M$ subject to $H$, i.e. a compact exhaustion $\{K_i\}$ such that $H([0,1]\times K_i) \subset \interior(K_{i+1})$. Consider the following diagram in $\OpenMan_d$
    \begin{center}
        \begin{tikzcd}
            K_1 \arrow{r}{gf} & K_2 \arrow{r}{gf} & K_3 \arrow{r}{gf} & K_4 \arrow{r}{gf} & \cdots
            \\
            K_1 \arrow{u}{\id} \arrow[hook]{r} & K_2 \arrow[hook]{r} \arrow{u}{\id} & K_3 \arrow[hook]{r} \arrow{u}{\id} & K_4 \arrow[hook]{r} \arrow{u}{\id} & \dots
        \end{tikzcd}
    \end{center}
    here, the horizontal morphisms denoted by $\hookrightarrow$ refer to the canonical inclusions $K_i \subset K_{i+1}$; furthermore, the square only commutes up to the chosen isotopy $H$. It is readily checked that applying colimit in $\OpenMan_d$ to the above diagram yields an embedding $M \hookrightarrow W_f$ isotopic to $\iota_0$. We now inductively apply the isotopy extension theorem \cite[Theorem 1.4, p.180]{Hirsch}, i.e. one square at a time, and hence we see that the above diagram, as an object in $\Fun(\mathbb{N},\Fun([1],\OpenMan_d))$, is equivalent to a strict diagram
    \begin{center}
        \begin{tikzcd}
            K_1 \arrow{r}{gf} & K_2 \arrow{r}{gf} & K_3 \arrow{r}{gf} & K_4 \arrow{r}{gf} & \cdots
            \\
            K_1 \arrow{u}{\varphi_1} \arrow[hook]{r} & K_2 \arrow[hook]{r} \arrow{u}{\varphi_2} & K_3 \arrow[hook]{r} \arrow{u}{\varphi_3} & K_4 \arrow[hook]{r} \arrow{u}{\varphi_4} & \dots
        \end{tikzcd}
    \end{center}
    where all the $\varphi_i$ are diffeomorphisms. As this diagram is strict, it induces on the colimits an actual diffeomorphism. Furthermore, $\iota_0$ is isotopic to this diffeomorphism, as the objects in the diagrams themselves are equivalent in $\Fun(\mathbb{N},\Fun([1],\OpenMan_d))$ prior to applying colimits. The fact that $\iota_0$ is isotopic to a diffeomorphism, together with diagram (\ref{diagram:emb-v-diff}), imply item (i).
    
    Regarding (ii), we begin with the following observation: given two group objects in spaces $G, H$, and a group map $f \colon G \to H$, $f$ is an equivalence if (a) it induces a $\pi_0$-surjection, and (b) the fibre of the identity is contractible. Let $W \in \OpenMan_d$, and fix a compact exhaustion $\{K_i\}_{i\in\mathbb{N}}$ of $W$. Without loss of generality, we can assume $W$ is noncompact, as otherwise the statement is obvious. We verify (a) and (b) of the above, for the map $\Diff(W) \to \Emb^\simeq(W,W)$. Note that (a) follows from (i), by placing $M=N=W$. Therefore, we show (b).
    
    \begin{itemize}
        \item \underline{\textit{Contractibility of the fibre of the identity embedding}}
    \end{itemize}
    Choose a compact exhaustion $\{K_i\}_{i \in \mathbb{N}}$ of $W$, and recall that $W \simeq \colim_{i\in\mathbb{N}}K_i$ in $\OpenMan_d$; thus, the map $\Diff(W) \to \Emb(W,W)\simeq \lim_i \Emb(K_i,W)$ is equivalent to the data of a map, for each $i \in \mathbb{N}$, of the form $\varphi_i\colon\Diff(W) \to \Emb(K_i,W)$, together with homotopies
    \begin{center}
        \begin{tikzcd}
            \Diff(W) \arrow[equal]{d} \arrow{r}{\varphi_{i+1}} & \Emb(K_{i+1},W) \arrow{d}{\res}
            \\
            \Diff(W) \arrow{r}[swap]{\varphi_i} & \Emb(K_i,W)
        \end{tikzcd}
    \end{center}
    where $\res$ is the restriction map along the fixed embeddings $K_i \subset K_{i+1}$. As $K_i$ is compact, parametrized isotopy extension holds, and one can identify the fibre of the map $\varphi_i$ at some fixed embedding $K_i \hookrightarrow W$ as
    $$\Diff_\partial(W\setminus\interior(K_i)) \to \Diff(W) \to \Emb(K_i,W)$$
    Fixing the identity basepoint in $\Emb(W,W)$ we obtain a fibre sequence
    $$\lim_{i\in\mathbb{N}}\Diff_\partial(W\setminus\interior(K_i)) \to \Diff(W) \to \Emb (W,W)_\id$$
    where the limit is obtained along the diagram $\{\Diff_\partial(W\setminus\interior(K_{i+1})) \to \Diff_\partial(W\setminus\interior(K_i))\}$ given by extension by identity. We claim that this (homotopy) limit is contractible.

    We have the following topological model for the homotopy limit
    $$\holim_i \Diff_\partial(W\setminus\interior(K_i)) \coloneqq \{(\varphi_i,\gamma_i)\in \prod_{i\in\mathbb{N}} \left(\Diff_\partial(W\setminus\interior(K_i))\right)\times\left(\Diff_\partial(W\setminus\interior(K_i))^I\right)\vert\,\gamma_i(0)=\varphi_i,\,\gamma_i(1)=\varphi_{i+1}\}$$
    where
    \begin{itemize}
        \item Each $\Diff_\partial(W\setminus\interior(K_i)$ is endowed with the weak Whitney $C^\infty$-topology,
        \item Each $\Diff_\partial(W\setminus\interior(K_i))^I$ is endowed with the compact-open topology,
        \item The infinite product is endowed with the product topology.
    \end{itemize}
    With respect to these topologies, we observe that the extension by identity morphism $\Diff_\partial(W\setminus\interior(K_{i+1}))\to \Diff_\partial(W\setminus\interior(K_i))$ is a homeomorphism onto its image, which allows us to reinterpret the whole tower as a sequence of subspace inclusions. For ease of notation, we impose the following:
    \begin{align*}
        &G_i \coloneqq \Diff_\partial(W\setminus\interior(K_i));
        \\
        &\Gamma_i \coloneqq \Diff_\partial(W\setminus \interior(K_i))^I;
        \\
        &X \coloneqq \holim_i\Diff_\partial(W\setminus\interior(K_i))\subset \prod_i(G_i \times \Gamma_i)
    \end{align*}
    We begin by writing a shift map. Define, for each $i \in \mathbb{N}$, maps
    \begin{align*}
        &\alpha_i \colon  \prod_{i\in\mathbb{N}}( G_i \times \Gamma_i) \to \prod_iG_i \xrightarrow[]{\mathrm{proj}_{i+1}}G_{i+1}\hookrightarrow G_i;
        \\
        &\beta_i \colon \prod_{i\in \mathbb{N}}(G_i \times \Gamma_i) \to \prod_i \Gamma_i \xrightarrow[]{\mathrm{proj}_{i+1}}\Gamma_{i+1} \hookrightarrow \Gamma_i
    \end{align*}
    and define $\sigma \coloneqq \prod_i(\alpha_i \times \beta_i)\colon \prod_i(G_i \times \Gamma_i) \to \prod_i(G_i \times \Gamma_i)$. By a routine check, one sees that $\sigma$ restricts to a map $\sigma \colon X \to X$. We first show that this map is homotopic to the identity. Define the map
    $$H\colon \prod_i(G_i \times \Gamma_i)\times I \to \prod_i(G_i \times \Gamma_i)$$
    by $((\varphi_i,\gamma_i),t) \mapsto (\gamma_i(t),\gamma_i \boxtimes\gamma_{i+1})$, where $\gamma_i\boxtimes\gamma_{i+1}$ denotes the image of the function $\prod_i\Gamma_i \times I \to \prod_i\Gamma_i$, sending the $i$-th factor to the function
    $$\gamma_i\boxtimes\gamma_{i+1}(s,t)=\begin{cases}
\gamma_i(s+t) & \text{if } t \leq 1-s, \\
\gamma_{i+1}(t-1+s) & \text{if } t \geq 1-s
\end{cases}$$
Note that for all $t$, we have that $\gamma_i\boxtimes\gamma_{i+1}(0,t)=\gamma_i(t)$, and $\gamma_i\boxtimes\gamma_{i+1}(1,t)=\gamma_{i+1}(t)$; therefore, the map $H$ restricts to a map $H\colon X \times I \to X$, starting at the identity, and ending at the shift map, implying that $\sigma$ is homotopic to the identity.

We now aim to reiterate this shift map. To this end, we consider the interval $[0,\infty)$, written as a union $\bigcup_{n\in \mathbb{N}\cup\{0\}}[n,n+1]$. We define a map $\Psi\colon X\times [0,\infty)\to X$, by imposing the function $\sigma^n\circ H$ on the interval $[n,n+1]$. In other words, we start with $H$ on $[0,1]$, $\sigma\circ H$ on $[1,2]$, $\sigma^2\circ H$ on $[2,3]$, and so on. It is clear that on the overlaps, the maps agree, and hence glue to a continuous function. Identifying $[0,\infty)$ with $[0,1)$, we have thus decomposed $[0,1)$ as a union of countably many adjacent intervals, which we denote $I_n$, and on each $I_n$, we write the function $\sigma^n\circ H$. We now define the function, which is a priori not continuous, denoted $\widetilde{\Psi}\colon X \times I \to X$, by defining it to be $\Psi$ on $X \times [0,\infty)$, and to be the constant map $(\id,\id)$ on $X\times\{1\}$. It remains to show that this map is continuous; once this is shown, this defines a retraction $X \simeq \ast$.

Let $U$ be an open neighborhood of $(\id,\id)$ in $X$; we may assume that $U = U_1 \times U_2$, where $U_1$ is open in $\prod_i G_i$ and $U_2$ is open in $\prod_i \Gamma_i$. By the definition of the product topology, $U_1$ constrains finitely many coordinates in $\prod_i G_i$, say those with coordinate $\leq j$ for some $j \in \mathbb{N}$. A basic open neighborhood of $G_i =\Diff_\partial (W \setminus\interior(K_i))$ is specified by a compact set $C \subset W\setminus\interior(K_i)$, an $\varepsilon>0$, such that all derivatives are $\varepsilon$-distant from the identity on $C$. Now, as $C$ is compact in $W\setminus\interior(K_i)$ and hence in $W$, it is fully contained in some $K_I$ for large $I$; in particular, it follows that all elements in $G_j=\Diff_\partial(W\setminus\interior(K_j))$ for $j\geq I$ are also contained in this open neighborhood, as they are all given by diffeomorphisms that are identity on $K_j$ which contains $C$. As there are finitely many coordinates with restraints, it follows that there is some $N$ large enough so that all $G_j$ are contained in the intersection of all the coordinates of $U_1$, for $j \geq N$. The same argument applied on $U_2$, after spelling out the definition of the compact-open topology, yields $N'$ so that all $\Gamma_j$ are contained in the intersection of all the projections of $U_2$. We thus let $n =\max\{N,N'\}$; then, it follows that $\widetilde{\Psi}(X \times ([0,1] \setminus I_N^\circ))\subset U$, implying continuity.
\end{proof}
As an immediate corollary, we obtain the following description of the groupoid cores of the categories $\CpctMan_d$ and $\OpenMan_d$.
\begin{cor}
    For all $d\in \mathbb{N}$, we have the following equivalences
    \begin{enumerate}
        \item $\left(\CpctMan_d\right)^\simeq \simeq \bigsqcup_{[W]}\BDiff(W)$, where $\{[W]\}$ is the set of diffeomorphism classes of smooth $d$-dimensional manifolds, which are diffeomorphic to the interior of a compact manifold.
        \item $\left(\OpenMan_d\right)^\simeq \simeq \bigsqcup_{[W]}\BDiff(W)$, where $\{[W]\}$ is the set of diffeomorphism classes of smooth $d$-dimensional manifolds with empty boundary.
    \end{enumerate}
\end{cor}

\subsection{Proof of \cref{thm:openMfldFreeColim}}
\subsubsection{The functor $\interior_!$}
By Lemma \ref{lem:OpenMan-sequential-colimits}, and by the universal property of sequential colimit completions, the functor $\interior\colon\CpctMan_d \hookrightarrow\OpenMan_d$ admits a unique, sequential colimit preserving extension $\interior_!\colon \PSh^\mathbb{N}(\CpctMan_d) \to \OpenMan_d$. We show that this functor is an equivalence of categories; for this, we employ Lemma \ref{lem:conditions}. The following is the missing condition, namely condition (ii) in Lemma \ref{lem:conditions}.

\begin{lem}\label{lem:Compact-manifolds-are-compact-objects}
    Let $W_\bullet\colon \mathbb{N} \to \OpenMan_d $ be a sequence of manifolds in $\OpenMan_d$, and denote its colimit in $\OpenMan_d $ by $W$. For $M\in \CpctMan_d $, the canonical map
    $$\colim_i \Emb (M,W_i) \to \Emb (M,W)$$
    is an equivalence.
\end{lem}
\begin{proof}
    We use point-set models to show this statement: recall that $\Emb (M,W)$ is equivalent to the homotopy type of the topological space $\underline{\Emb }(M,W)$, endowed with the (weak) $C^\infty$-topology \cite[Chapter 2]{Hirsch}. We will use the notation $\underline{\Emb }$ to denote the topological space of embeddings. We list some facts, which are immediate consequences of the definition of the corresponding topology on the spaces of embeddings:
    \begin{itemize}
        \item We have inclusions of subsets $\underline{\Emb }(M,W_i) \subset \underline{\Emb }(M,W_{i+1})$ induced by the embeddings $W_i \hookrightarrow W_{i+1}$, and as $M \in \CpctMan_d$, we have
        $$\underline{\Emb }(M,W)= \bigcup_{i \in \mathbb{N}}\underline{\Emb }(M,W_i)$$
        \item $\underline{\Emb }(M,W_i)$, as a topological space endowed with the subspace topology inherited from $\underline{\Emb }(M,W_{i+1})$, is homeomorphic to the original topological space $\underline{\Emb }(M,W_i)$;
        \item The map $\underline{\Emb }(M,W_i) \hookrightarrow \underline{\Emb }(M,W_{i+1})$ is an open map;
        \item The topological space $\underline{\Emb }(M,W)$ is homeomorphic to the topological space
        $$\bigcup_{i \in \mathbb{N}}\underline{\Emb }(M,W_i)$$
        endowed with the union topology.
    \end{itemize}
    It is a general fact that given a topological space $X$ with an increasing sequence of open subsets $\{\mathcal{U}_i\}_{i \in \mathbb{N}}$, the topological space $\bigcup_i \mathcal{U}_i$ endowed with the union topology is a point-set model for the homotopy colimit $\hocolim_i \mathcal{U}_i$, in the sense that the map
    $$\hocolim_i \mathcal{U}_i \to \bigcup_i \mathcal{U}_i$$
    is a weak equivalence.
\end{proof}
We are now ready to prove the first part of the main theorem.
\begin{proof}[Proof of Theorem \ref{thm:openMfldFreeColim}, part 1]
    We show that the functor $\interior_!\colon \PSh^\mathbb{N}(\CpctMan_d) \to \OpenMan_d$ is an equivalence. We follow Lemma \ref{lem:conditions}. Condition (i) is the fact that $\interior\colon \CpctMan_d \to \OpenMan_d$ is fully faithful, which was established in Lemma \ref{lem:interior-functor-fully-faithful}. Condition (ii) is established in Lemma \ref{lem:Compact-manifolds-are-compact-objects}, and condition (iii) is the fact that any manifold admits a compact exhaustion.
\end{proof}
\subsubsection{The inverse of $\interior_!$ as restricted Yoneda}
We now show that the inverse of the above functor admits a description as restricted Yoneda. By Lemma \ref{lem:right-adjoint-colim-extension}, it is sufficient to show that the restriction functor $\interior^\ast \colon \PSh(\OpenMan_d) \to \PSh(\CpctMan_d)$ restricts to a functor $\interior^\ast\colon \PSh^\mathbb{N}(\OpenMan_d) \to \PSh^\mathbb{N}(\CpctMan_d)$; once done, Lemma \ref{lem:right-adjoint-colim-extension} ensures that the functor
$$\Exh\colon \OpenMan_d \hookrightarrow \PSh^\mathbb{N}(\OpenMan_d) \xrightarrow{\interior^\ast} \PSh^\mathbb{N}(\CpctMan_d)$$
is indeed the right adjoint to $\interior_!$, and hence its inverse.
\begin{lem}\label{lem:colimits}
    Let $W \in \OpenMan_d $, and let $\{W_i\}_{i \in \mathbb{N}}$ be an increasing collection of open submanifolds of $W$, such that $W = \bigcup_i W_i$, i.e. $W\simeq \colim_i W_i \in \OpenMan_d $. The natural morphisms $\{h^\open_{W_i}\} \to h^\open_W$ assemble into an equivalence
    $$\colim_i \left(\interior^\ast h^\open_{W_i}\right) \to \interior^\ast h^\open_W$$
    in $\PSh(\CpctMan_d )$. Here, $h^\open$ denotes the Yoneda embedding $\ManInf\hookrightarrow \PSh^\mathbb{N}(\ManInf_d)$.
\end{lem}
\begin{proof}
    We expand slightly on the constructions of the maps in the statement above. By definition, each $W_i$ is equipped with a morphism $W_i \hookrightarrow W$ in $\OpenMan_d $, and these morphisms are compatible with the maps $W_i \hookrightarrow W_{i+1}$; in particular, we obtain a collection of compatible maps $h^\open_{W_i} \to h^\open_W$ in $\PSh(\OpenMan_d )$. Applying $\interior^\ast$ to the collection of morphisms $\{h^\open_{W_i} \to h^\open_W\}$, we obtain a collection of morphisms $\interior^\ast h^\open_{W_i} \to \interior^\ast h^\open_W$, compatible with all the maps $\{\interior^\ast h^\open_{W_i} \to h^\open_{W_{i+1}}\}$. In particular, these morphisms assemble into a map
    $$\colim_i \left(\interior^\ast h^\open_{W_i}\right) \to \interior^\ast h^\open_W$$
    which we would like to show is an equivalence. As this is a map of presheaves, this will be the case if this map induces an equivalence at all evaluations at $M \in \CpctMan_d $. Writing out what these presheaves evaluate to, we obtain a map
    \begin{equation}\label{eq:colim}
        \colim_i \Emb (M, W_i) \to \Emb (M,W)
    \end{equation}
    which is an equivalence by Lemma \ref{lem:Compact-manifolds-are-compact-objects}.
\end{proof}
With the above lemma at hand, we are one step closer to defining the restricted Yoneda functor; the following lemma is the final input for this construction.
\begin{lem}\label{lem:restriction-functor-exists}
    The functor $\interior^\ast \colon \PSh(\OpenMan_d ) \to \PSh(\CpctMan_d )$ restricts to a functor as suggested in the following commutative diagram
    \begin{center}
        \begin{tikzcd}
            \PSh^\mathbb{N}(\OpenMan_d ) \arrow[hook]{d} \arrow[dashed]{r} & \PSh^\mathbb{N}(\CpctMan_d ) \arrow[hook]{d}
            \\
            \PSh(\OpenMan_d ) \arrow{r}[swap]{\interior^\ast} & \PSh(\CpctMan_d )
        \end{tikzcd}
    \end{center}
\end{lem}
\begin{proof}
    We begin by observing that as $\PSh^\mathbb{N}(\CpctMan_d ) \hookrightarrow \PSh(\CpctMan_d )$ is a full subcategory inclusion, it is a property that the functor 
    $$\PSh^\mathbb{N}(\OpenMan_d ) \hookrightarrow \PSh(\OpenMan_d ) \xrightarrow[]{\interior^\ast} \PSh(\CpctMan_d )$$
    lands in $\PSh^\mathbb{N}(\CpctMan_d )$. By Lemma \ref{lem:naive-colim}, it suffices to show the following: given an object $X\in\PSh(\OpenMan_d )$, which is a sequential colimit of representables $\{h^\open_{W_i}\}$, then $\interior^\ast X$ is equivalent to a sequential colimit of representables in $\PSh(\CpctMan_d )$. We therefore fix such an $X$. As $\interior^\ast$ preserves colimits, we obtain
    $$\interior^\ast X \simeq \interior^\ast \colim_{i\in\mathbb{N}}h^\open_{W_i} \simeq \colim_{i\in \mathbb{N}}\interior^\ast h^\open_{W_i}$$
    For each $i \in \mathbb{N}$, we choose an exhaustion $\{K_{i,j}\}_{j \in \mathbb{N}}$ of $W_i$ by compact, codimension 0 submanifolds. Then, by Lemma \ref{lem:colimits}, we obtain an equivalence
    $$\colim_{j \in \mathbb{N}}\left(\interior^\ast h^\open_{K_{i,j}}\right) \simeq \interior^\ast h^\open_{W_i}$$
    We observe that for each $j$, $K_{i,j}$ was assumed to be compact, and hence
    $$\interior^\ast h^\open_{K_{i,j}} \simeq h^{\mathrm{cpt}}_{K_{i,j}}$$
    where $h^{\mathrm{cpt}}$ denotes the Yoneda embedding for $\CpctMan_d$. Combining all the above, we obtain the following chain of equivalences:
    \begin{align*}
        \interior^\ast X & \simeq \colim_{i \in \mathbb{N}}\left( \interior^\ast h^\open_{W_i} \right)\simeq \colim_{i \in \mathbb{N}} \left(\colim_{j \in \mathbb{N}}h^\cpct_{K_{i,j}}\right)
    \end{align*}
    By the same argument as in the proof of Lemma \ref{lem:naive-colim}, it follows that $\interior^\ast X$ can be represented as a sequential colimit of representables in $\PSh(\CpctMan_d )$.
\end{proof}
We are now ready to define the restricted Yoneda functor.
\begin{dfn}\label{dfn:Exh}
    Define $\Exh$ as the following composition
    \begin{center}
        \begin{tikzcd}
            \Exh\colon\OpenMan_d  \arrow[hook]{r} & \PSh^\mathbb{N}(\OpenMan_d ) \arrow{r}{\interior^\ast} & \PSh^\mathbb{N}(\CpctMan_d )
        \end{tikzcd}
    \end{center}
    where $\interior^\ast$ exists by virtue of Lemma \ref{lem:restriction-functor-exists}.
\end{dfn}

The notation $\mathrm{Exh}$ refers to the rather informal observation that this functor mirrors the construction of a compact exhaustion of a manifold by compact, codimension 0 submanifolds.
We are now ready to prove the second part of Theorem \ref{thm:openMfldFreeColim}, namely that $\Exh$ is the inverse of the functor $\interior_!$.

\begin{proof}[Proof of Theorem \ref{thm:openMfldFreeColim}, part 2]
   By Lemma \ref{lem:right-adjoint-colim-extension}, we know that $\interior_!$ is left adjoint to $\Exh$; we have shown that $\interior_!$ is an equivalence, and therefore $\Exh$ is indeed its inverse.
\end{proof}

As a corollary of Theorem \ref{thm:openMfldFreeColim}, one obtains the following universal property. This is merely a write-up for the universal property of $\PSh^\mathbb{N}(-)$ in the setting of categories of manifolds.
\begin{cor}\label{cor:universal-property-manifolds}
    Let $\catsingle{C}$ be an $\infty$-category that admits sequential colimits. Then, restriction along $\interior\colon\CpctMan_d  \hookrightarrow\OpenMan_d $ induces an equivalence
    $$\Fun^\mathbb{N}(\OpenMan_d ,\catsingle{C}) \to \Fun(\CpctMan_d ,\catsingle{C})$$
    where $\Fun^\mathbb{N}$ denotes the category of functors preserving sequential colimits.
\end{cor}

\section{Bordism-Ends}
\subsection{A pullback square}\label{section:endbordisms}
For notation, we let $\BordInf_d $ denote the $(\infty,1)$-category of closed $(d-1)$-dimensional manifolds, and bordisms between them. We follow the construction of \cite{discstructurespace,AdelaJan}; as these references are extensive, we do not recall the lengthy definition here. One way to relate $\BordInf_d $ to the categories $\CpctMan_d $ studied above can be extracted from \cite[Corollary A.2]{AdelaJan}: there is an equivalence of categories
$$\Phi\colon\left(\BordInf_d \right)_{\emptyset/} \simeq \CpctMan_d $$
which, on objects, views a compact $d$-manifold $M$ as a bordism $\emptyset \rightsquigarrow\partial M$, and on morphisms, sends an embedding $e\colon M \hookrightarrow N$ to the bordism $N\setminus\interior(e(M))\colon\partial M \rightsquigarrow \partial N$, making the following triangle commute
\begin{center}
    \begin{tikzcd}
        & \emptyset \arrow{dl}[swap]{M} \arrow{dr}{N}
        \\
        \partial M \arrow{rr}[swap]{N\setminus \interior(e(M))} && \partial N
    \end{tikzcd}
\end{center}
With this equivalence at hand, we define the following 
\begin{dfn}\label{def:boundary-bordismends-functors}
We define the following two functors.
\begin{itemize}
    \item Let $\partial \colon \CpctMan_d \to \BordInf_d$ be the functor obtained as the composition
    $$\partial\colon \CpctMan_d \xrightarrow[\simeq]{\Phi} (\BordInf_d)_{\emptyset/} \xrightarrow[]{\mathrm{fgt}} \BordInf_d$$
    \item Let $\bordEnd\colon \OpenMan_d \to \PSh^\mathbb{N}(\BordInf_d)$ be the functor obtained as the composition
    $$\bordEnd\colon\OpenMan_d \xrightarrow[\simeq]{\PSh^\mathbb{N}(\Phi)} \PSh^\mathbb{N}(\left(\BordInf_d\right)_{\emptyset/}) \xrightarrow{\PSh^\mathbb{N}(\mathrm{fgt})} \PSh^\mathbb{N}(\BordInf_d)$$
\end{itemize}
\end{dfn}
For a manifold $W$, we call the object $\bordEnd(W)$ \emph{the bordism-end of $W$}. If a compact exhaustion $\{K_i\}_{i\in \mathbb{N}}$ is chosen, a routine check shows that $\bordEnd(W)$ is equivalent, in $\PSh^\mathbb{N}(\BordInf_d)$, to the system $\{\partial K_i\}$ with structure maps $K_i \setminus\interior(K_{i-1})$.

\begin{thm}\label{thm:Ind-End}
    For all $d\in \mathbb{N}$, the following commutative square is a pullback square of $\infty$-categories
    \begin{center}
        \begin{tikzcd}
            \CpctMan_d \arrow[hook]{r}{\interior} \arrow{d}[swap]{\partial} & \OpenMan_d \arrow{d}{\bordEnd}
            \\
            \BordInf_d \arrow[hook]{r} & \PSh^\mathbb{N}(\BordInf_d)
        \end{tikzcd}
    \end{center}
    
\end{thm}\label{thm:pullback-square}
\begin{proof}
    The above square is, by definition, the outer square of the following commutative diagram
    \begin{center}
        \begin{tikzcd}
            \CpctMan_d  \arrow[hook]{r} \arrow{d}{\simeq}[swap]{\Phi} & \OpenMan_d  \simeq \PSh^\mathbb{N}(\CpctMan_d ) \arrow{d}{\PSh^\mathbb{N}(\Phi)}[swap]{\simeq}
            \\
            \left(\BordInf_d \right)_{\emptyset/} \arrow[hook]{r} \arrow{d}[swap]{\mathrm{fgt}} & \PSh^\mathbb{N}\left(\left(\BordInf_d \right)_{\emptyset/}\right) \arrow{d}{\PSh^\mathbb{N}(\mathrm{fgt})}
            \\
            \BordInf_d  \arrow[hook]{r} & \PSh^\mathbb{N}\left(\BordInf_d \right)
        \end{tikzcd}
    \end{center}
    We show that each of the squares are pullbacks. Observe that the first square is clearly a pullback square: the functor $\Phi$ is an equivalence of categories, and since we saw that $\PSh^\mathbb{N}$ defines an endofunctor on the $\infty$-category of $\infty$-categories, in particular $\PSh^\mathbb{N}(\Phi)$ is also an equivalence.

    Regarding the lower square, we make reference to Lemma \ref{lem:slice} to show that the square is equivalent to the following square
    \begin{center}
        \begin{tikzcd}
            \left(\BordInf_d \right)_{\emptyset/} \arrow[hook]{r} \arrow{d}[swap]{\mathrm{fgt}} & \PSh^\mathbb{N}\left(\left(\BordInf_d \right)\right)_{\emptyset/} \arrow{d}{\mathrm{fgt}}
            \\
            \BordInf_d  \arrow[hook]{r} & \PSh^\mathbb{N}\left(\BordInf_d \right)
        \end{tikzcd}
    \end{center}
    which clearly is a pullback square, as the horizontal functors are fully faithful.
\end{proof}
For $W$ a manifold the above theorem is equivalent to the following, more explicit, statement: $W$ is diffeomorphic to the interior of a compact manifold with boundary if and only if the object $\bordEnd(W)$ is equivalent to a representable object in $\PSh^\mathbb{N}(\BordInf_d)$; that representable object (which is now a closed $(d-1)$-manifold in $\BordInf_d$) is the boundary of that compact manifold. This boundary is only well-defined up to equivalences in $\BordInf_d$, i.e. invertible bordisms, as expected.

\subsection{A fibre sequence}
In this section, we implement the bordism-theoretic perspective of §\ref{section:endbordisms} to obtain a fibre sequence relating moduli spaces of manifolds with their compactly supported analogues.
\subsubsection{A mapping space description of the fibre}
We begin with a general observation. Let $\catsingle{C}$ be an $\infty$-category, and fix an object $c\in\catsingle{C}$. Then, we obtain a pullback square in $\CatInf$
\begin{center}
    \begin{tikzcd}
        \Map_\catsingle{C}(c,d) \arrow{r}\arrow{d} & \catsingle{C}_{c/}\arrow{d}{\forget}
        \\
        \ast \arrow{r} & \catsingle{C}
    \end{tikzcd}
\end{center}
where $\Map_{\catsingle{C}}(c,d)$ is viewed as an $\infty$-groupoid. We now apply the above general formalism to the setting of Theorem \ref{thm:Ind-End}. By Theorem \ref{thm:Ind-End} and Lemma \ref{lem:slice}, we obtain an equivalence $\OpenMan_d  \simeq \PSh^\mathbb{N}(\BordInf_d )_{\emptyset/}$, where $\emptyset$ is viewed as an object in $\PSh^\mathbb{N}(\BordInf_d )$ via the Yoneda embedding. In Theorem \ref{thm:Ind-End}, a functor $\OpenMan_d  \to \PSh^\mathbb{N}(\BordInf_d )$ played a central role, and it follows readily that this functor, under the identification $\OpenMan_d  \simeq \PSh^\mathbb{N}(\BordInf_d )_{\emptyset/}$, corresponds to the forgetful functor
$$\PSh^\mathbb{N}(\BordInf_d )_{\emptyset/} \to \PSh^\mathbb{N}(\BordInf_d )$$
We thus obtain a pullback square in $\CatInf$
\begin{equation}
    \begin{tikzcd}
        \Map_{\PSh^\mathbb{N}(\BordInf_d )}(\emptyset,Z) \arrow{r} \arrow{d} & \OpenMan_d  \arrow{d}{\bordEnd}
        \\
        \ast \arrow{r}[swap]{\{Z\}} & \PSh^\mathbb{N}(\BordInf_d )
    \end{tikzcd}
\end{equation}
Applying groupoid cores, we obtain a pullback
\begin{equation}\label{diagram:mapping-space-fibre}
    \begin{tikzcd}
        \Map_{\PSh^\mathbb{N}(\BordInf_d )}(\emptyset,Z) \arrow{r} \arrow{d} & \left(\OpenMan_d \right)^\simeq \arrow{d}{\bordEnd^\simeq}
        \\
        \ast \arrow{r}[swap]{\{Z\}} & \left(\PSh^\mathbb{N}(\BordInf_d )\right)^\simeq
    \end{tikzcd}
\end{equation}
\subsubsection{Compactly supported diffeomorphisms}
We now investigate the mapping space fibre in diagram (\ref{diagram:mapping-space-fibre}); the claim is that this will informally model a "compactly supported moduli space of manifolds". The goal here is to make this statement more precise. To that end, we begin with the following definition.
\begin{dfn}
    Let $W\in \OpenMan_d$, and choose a compact exhaustion $\{K_i\}_{i\in \mathbb{N}}$. Define the group object in $\catsingle{S}$, $\Diff_c(W)$, as the colimit in $\catsingle{S}$ of $\{\Diff_\partial(K_i)\}$, where the structure maps are induced by extensions by identity.
\end{dfn}
We first observe that the above is well defined, i.e. is independent of the choice of compact exhaustions. If $\{K_i\}$ and $\{K_i'\}$ are two choices of compact exhaustions, one can find subsequences $\{K_\varphi(i)\}$ and $\{K_{\psi(i)}'\}$ for $\varphi,\psi\colon \mathbb{N} \to \mathbb{N}$ increasing, unbounded functions, such that for all $i$, $K_{\varphi(i)}\subset K_{\psi(i)}' \subset K_{\varphi(i+1)}$. Thus, rechoosing a compact exhaustion interleaving $K_i$ and $K_j'$ as above, it follows that both colimits are equivalent, both equivalent to the colimit over the new tower hence defined.
\begin{rem}
    We justify why we call this colimit the "space of compactly supported diffeomorphisms".
\begin{enumerate}
    \item One could consider the colimit of the tower $\{\Diff_\partial (K_i)\}$, with the Whitney $C^\infty$-topology, along the extension by identity maps. The extension by identity maps are cofibrations: $\Diff_\partial(K_i)$ is closed in $\Diff_\partial(K_{i+1})$, and they are both ANR's as a corollary of \cite[§4, Proposition 1 (2)]{diffcompactsupport} and \cite[Corollary to Theorem 5]{Palais}, and for ANR's, closed embeddings are cofibrations \cite[Theorem 4.2.15]{ANR}. Thus, the colimit in topological spaces is homotopy equivalent to the homotopy colimit.
    \item One could define $\Diff_c(W)$ as a subspace of $\Diff(W)$, endowed with the strong Whitney topology. A straightforward argument shows that one obtains a continuous map
    $$\colim \Diff_\partial (K_i) \to \Diff_c(W)$$
    where the left-hand-side is the colimit topology from (i), and one observes that the above map induces a weak homotopy equivalence. This rests on the following observations: (a) the map $\Diff_\partial (K_i) \to \Diff(W)$, given by extension by identity, is a homeomorphism onto its image; (b) any compact subset $K \subset \Diff_c(W)$ is fully contained in $\Diff_\partial (K_i)$ for large enough $i$. Both the above statements follow readily from the definition of the strong Whitney topology;
    \item The colimit topology can also be shown to be weakly homotopy equivalent to another topology on $\Diff_c(W)$, namely to the very strong Whitney topology, (refered to the Whitney topology in \cite{diffcompactsupport}). In other words, this space appears as the underlying homotopy type of a topology studied in literature on spaces of compactly supported diffeomorphisms.
\end{enumerate}
\end{rem}

We can now give a more concrete description of the space $\Map_{\PSh^\mathbb{N}(\BordInf_d )}(\emptyset,\bordEnd (W))$. Indeed, we have the following equivalences.
\begin{align*}
    \Map_{\PSh^\mathbb{N}(\BordInf_d )}(\emptyset,\bordEnd (W)) &\simeq \Map_{\PSh^\mathbb{N}(\BordInf_d )}(\emptyset,\colim_{i\in\mathbb{N}}\partial K_i)
    \\
    &\simeq\colim_{i\in\mathbb{N}}\Map_{\BordInf_d }(\emptyset,\partial K_i)
    \\
    &\simeq \colim_{i\in\mathbb{N}}\left(\bigsqcup_{[W_i,\partial W_i\cong\partial K_i]}\BDiff_\partial(W_i)\right)
\end{align*}
where $[W_i]$ in the last line denotes the set of diffeomorphism classes rel $\partial$ of compact manifolds with a fixed identification of their boundaries with $\partial K_i$. As $\pi_0$ commutes with sequential colimits, we obtain that
$$\pi_0 \Map_{\PSh^\mathbb{N}(\BordInf_d)}(\emptyset,Z)\cong \colim_{i\in\mathbb{N}}\pi_0\Map_{\BordInf_d}(\emptyset,\partial K_i)$$
We now describe each component. To that end, fix a map $\ast \to \pi_0 \Map_{\PSh^\mathbb{N}(\BordInf_d)}(\emptyset,Z) \cong \colim_i \pi_0 \Map_{\BordInf_d}(\emptyset,\partial K_i)$. In particular, this is equivalent to fixing a component of $\Map_{\BordInf_d}(\emptyset,\partial K_j)$ for large enough $j$, i.e. a class of a diffeomorphism rel. boundary of some compact manifold $W_j$ with boundary identified with $\partial K_j$. This component then has the following homotopy type:
$$\Map_{\PSh^\mathbb{N}(\BordInf_d)}(\emptyset,Z)_\ast \simeq \colim_{i\geq j}\BDiff_\partial(W_i)$$
where $W_i$ is obtained from $W_j$ by attaching the cobordisms $\partial K_j \rightsquigarrow \partial K_i$ described by $Z$. We have the equivalences
\begin{align*}
    \colim_{i\geq j} \BDiff_\partial (W_i) &\simeq \mathrm{B}(\colim_{i\geq j}\Diff_\partial (W_i)) \simeq \BDiff_c(W)
\end{align*}
where the first equivalence follows from the fact that $B$ is a left adjoint. Here, $W \coloneqq \colim_{i\geq j}W_i$ in $\OpenMan_d$.

Finally, we identify $\pi_0 \Map_{\PSh^\mathbb{N}(\BordInf_d)}(\emptyset,Z)$ more geometrically: recall that we have a pullback square
\begin{center}
    \begin{tikzcd}
        \Map_{\PSh^\mathbb{N}(\BordInf_d)}(\emptyset,Z) \arrow{r}\arrow{d} & \OpenMan_d \arrow{d}{\bordEnd}
        \\
        \ast \arrow{r}[swap]{\{Z\}} & \PSh^\mathbb{N}(\BordInf_d)
    \end{tikzcd}
\end{center}
Thus, $\pi_0\Map_{\PSh^\mathbb{N}}(\emptyset,Z)$ can be described as the set of diffeomorphism classes of manifolds $W \in \OpenMan_d$, with an identification $\bordEnd(W) \simeq Z$, up to diffeomorphisms rel. the $\bordEnd$ identification. Putting the above discussion together, we obtain the following
\begin{lem}
    Let $Z \in \PSh^\mathbb{N}(\BordInf_d)$. Then, we have an equivalence
    $$\Map_{\PSh^\mathbb{N}(\BordInf_d)}(\emptyset,Z) \simeq \bigsqcup_{[W,\bordEnd(W)\simeq Z]}\BDiff_c(W)$$
\end{lem}
\subsubsection{The fibre sequence}
We can now put the above discussion together to obtain the following theorem.
\begin{thm}\label{thm:compactly-supported-diff}
    The following is a pullback of spaces
    \begin{center}
        \begin{tikzcd}
            \bigsqcup_{[W,\bordEnd(W)\simeq Z]}\BDiff_c(W) \arrow{r}\arrow{d} & \bigsqcup_{[W]}\BDiff(W)\arrow{d}{\bordEnd^\simeq}\simeq \left(\OpenMan_d\right)^\simeq
            \\
            \ast \arrow{r}[swap]{\{Z\}} & \PSh^\mathbb{N}(\BordInf_d)^\simeq
        \end{tikzcd}
    \end{center}
    where $[W]$ ranges over diffeomorphism classes of manifolds with empty boundary, and where $[W,\bordEnd(W)\simeq Z]$ ranges over pairs of a manifold $W$ and an identification $\bordEnd(W)\simeq Z$, up to diffeomorphism rel. the bordism-end identification.
\end{thm}
\begin{rem}
    Restricting the functor $\bordEnd\colon \OpenMan_d \to \PSh^\mathbb{N}(\BordInf_d)$ to $\CpctMan_d$, we obtain the functor $\partial\colon \CpctMan_d \to \BordInf_d$ from Theorem \ref{thm:Ind-End}. The pullback of the span $\CpctMan_d \xrightarrow[]{\partial} \BordInf_d \xleftarrow[]{\{M\}}\ast$ is equivalent to the space $\Map_{\BordInf_d}(\emptyset,M)$, and by definition, this is equivalent to 
    $$\Map_{\BordInf_d}(\emptyset,M) \simeq \bigsqcup_{[W,\partial W\cong M]}\BDiff_\partial(W)$$
    where the components range over the set of compact, smooth $d$-manifolds with an identification $\partial \cong M$, up to diffeomorphisms rel. $\partial$ (note that this is the empty space if $M$ is not null-cobordant). In particular, one obtains a fibre sequence
    $$\bigsqcup_{[W,\partial W\cong M]}\BDiff_\partial (W) \to \bigsqcup_{[\interior(W)]} \BDiff(\interior(W)) \to \left(\BordInf_d\right)^\simeq_{M}$$
    where the base denotes the component of $M$ in the groupoid core of $\BordInf_d$, and where the components of the total space range over the diffeomorphism classes of interiors of compact, smooth $d$-manifolds. In the case where $d \geq 5$, $\left(\BordInf_d\right)^\simeq_M$ admits a more classical description: it is a delooping of the space of self h-cobordisms of $M$. This is the case for $d\geq 6$ by virtue of the s-cobordism theorem, and by \cite[Theorem 4]{Stallings}, \cite[Remark 6.19]{ManuelSanderInfinityOperadic} for $d=5$.
\end{rem}

\section{$\Pro$-space of ends}\label{section:pro-ends}
Given a manifold $W$, a more classical method to record the "end" of the manifold is as follows: a choice of a compact exhaustion $\{K_i\}_{i\in\mathbb{N}}$ gives rise to an inverse system of spaces, namely $\{W\setminus K_i\}$. Our goal in this section is to assemble this into a functor $\proEnd\colon\OpenMan_d \to \Cosp(\Pro(\catsingle{S}))$.
\subsection{The tail functor}\label{section:tail-functor}
We begin by constructing a functor
\begin{equation}\label{tail-functor}
    \tau\colon \PSh^\mathbb{N}(\Cosp(\catsingle{S})) \to \Cosp(\Pro(\catsingle{S}))
\end{equation}
This section follows \cite{Spans}, and we refer to this work for all the definitions and conventions. For $\catsingle{C}$ a category admitting pullbacks, we take $\mathrm{Span}(\catsingle{C})$ to be the $\infty$-category of spans in $\catsingle{C}$, informally described by the category where objects are objects of $\catsingle{C}$, and where a morphism $c \rightsquigarrow d$ is a span $c \leftarrow e \rightarrow d$ in $\catsingle{C}$; composition is obtained by pullbacks; for a more detailed definition, we refer to \emph{loc. cit.}. If $\catsingle{C}$ is an $\infty$-category admitting pushouts, we let $\Cosp(\catsingle{C})\coloneqq \mathrm{Span}(\catsingle{C}^\op)$. Informally, $\Cosp(\catsingle{C})$ is the category where objects are objects of $\catsingle{C}$, and where a morphism $c \rightsquigarrow d$ is a cospan $c \rightarrow e \leftarrow d$ in $\catsingle{C}$; composition is obtained by pushouts.

We now proceed with the actual construction of functor (\ref{tail-functor}). By Theorem \ref{thm:localisation}, our strategy is to first define a functor $\Fun(\mathbb{N},\catsingle{S}) \to \Cosp(\Pro(\catsingle{S}))$, which furthermore sends all the cofinal re-indexing morphisms $W$ in $\Fun(\mathbb{N},\catsingle{S})$ to equivalences; once this is done, the functor $\tau$ is provided by the universal property of localizations.
\begin{construction}
    We now construct a functor 
    $$\mathrm{tail}\colon \Fun(\mathbb{N},\Cosp(\catsingle{S})) \to \Cosp(\Pro(\catsingle{S}))$$
    By \cite[Theorem 2.18]{Spans}, there is an equivalence of $\infty$-categories
    $$\Fun(\mathbb{N},\Cosp(\catsingle{S})) \simeq \mathrm{Span}(\Fun_{\mathrm{AdTrip}}(\tw^r(\mathbb{N}),\catsingle{S}^\op))$$
    where the twisted arrow category $\tw^r(\mathbb{N})$ is endowed with its adequate triple structure as in \cite[Example 2.9]{Spans}, and where $\catsingle{S}^\op$ is endowed with the adequate triple structure $(\catsingle{S}^\op,\catsingle{S}^\op,\catsingle{S}^\op)$, i.e. all maps are both ingressive and egressive Example 2.9 (2) of \emph{loc. cit.}; furthermore, $\Fun_\mathrm{AdTrip}$ is itself endowed with the adequate triple structure as in Definition 2.16 of \emph{loc. cit.}, and $\mathrm{Span}$ is understood as the image of the functor $\mathrm{Span}\colon \mathrm{AdTrip} \to \CatInf$. The category of functors of adequate triples is by definition a full subcategory of the functor category $\Fun(\tw^r(\mathbb{N}),\catsingle{S}^\op)$, namely as the full subcategory on those functors that send ambigressive pullbacks to pullbacks, see the discussion following Definition 2.1 in \emph{loc. cit.} In particular, the fully faithful inclusion
    $$\inc\colon\Fun_{\mathrm{AdTrip}}(\tw^r(\mathbb{N}),\catsingle{S}^\op) \to \Fun(\tw^r(\mathbb{N}),\catsingle{S}^\op)$$
    can itself be viewed as a functor of adequate triples, where the source is given its adequate triple structure from earlier, and the target is endowed with the constant adequate triple structure, as $\Fun(\tw^r(\mathbb{N}),\catsingle{S}^\op)$ admits pullbacks (taken pointwise). As limits in the target are computed pointwise, it in particular sends ambigressive pullbacks to pullbacks. Our final observation is the following: if $\catsingle{C}$ is a category admitting pullbacks, then its usual span category $\mathrm{Span}(\catsingle{C})$ is defined as $\mathrm{Span}((\catsingle{C},\catsingle{C},\catsingle{C}))$ of the constant adequate triple. Combining the above, we obtain from the functor $\mathrm{Span}\colon \mathrm{AdTrip} \to \CatInf$, a functor
    $$\mathrm{Span}(\inc)\colon\mathrm{Span}(\Fun_\mathrm{AdTrip}(\tw^r(\mathbb{N}),\catsingle{S}^\op)) \to \mathrm{Span}(\Fun(\tw^r(\mathbb{N}),\catsingle{S}^\op))$$
    where the target is the usual span category of the functor category $\Fun(\tw^r(\mathbb{N}),\catsingle{S}^\op)$, where pullbacks are computed pointwise. We therefore now switch completely to the latter.

    Let $\mathbb{I}$ denote the poset of intervals $[i,j]$, with $i,j \in \mathbb{N}$, partially ordered by inclusion; that is, we have a morphism $[i,j] \to [i',j']$ precisely when $i'\leq i\leq j\leq j'$. By direct inspection, we have an equivalence $\tw^r(\mathbb{N}) \simeq \mathbb{I}^\op$; thus, $\Fun(\tw^r(\mathbb{N}),\catsingle{S}^\op) \simeq \Fun(\mathbb{I},\catsingle{S})^\op$. As a side-remark, we note that the full subcategory $\Fun_\mathrm{AdTrip}(\tw^r(\mathbb{N}),\catsingle{S}^\op)$ can be readily characterized as the following subcategory of $\Fun(\mathbb{I},\catsingle{S})^\op$: indeed, $X_{\bullet,\bullet}\colon\mathbb{I} \to \catsingle{S}$ is in $\Fun_\mathrm{AdTrip}(\tw^r(\mathbb{N},\catsingle{S}^\op))$ if and only if for all $j\geq i+2$, the natural comparison map $X_{i,j}\to X_{i,j-1}\times_{X_{i+1,j-1}} X_{i+1,j}$ is an equivalence. This can be readily seen from expanding the definitions.

    Let $s\colon\mathbb{I} \to \mathbb{N}^\op$ denote the functor sending $[i,j]$ to $i \in \mathbb{N}$. We then consider the left Kan extension functor $s_! \colon \Fun(\mathbb{I},\catsingle{S}) \to \Fun(\mathbb{N}^\op,\catsingle{S})$. Finally, we have a functor $\lim\colon\Fun(\mathbb{N}^\op,\catsingle{S}) \to \Pro(\catsingle{S})$, sending an inverse system of spaces to its limit in $\Pro(\catsingle{S})$.

    We now finally define $\mathrm{tail}$ as the following composite:
   
    \begin{center}
        \begin{tikzcd}
            \Fun(\mathbb{N},\Cosp(\catsingle{S})) \simeq \mathrm{Span}(\Fun_\mathrm{AdTrip}(\tw^r(\mathbb{N}),\catsingle{S}^\op)) \arrow{d}{\mathrm{Span}(\inc)}
            \\
            \mathrm{Span}(\Fun(\tw^r(\mathbb{N}),\catsingle{S}^\op))\simeq \mathrm{Span}(\Fun(\mathbb{I},\catsingle{S})^\op) = \Cosp(\Fun(\mathbb{I},\catsingle{S})) \arrow{d}{\Cosp(s_!)}
            \\
            \Cosp(\Fun(\mathbb{N}^\op,\catsingle{S})) \arrow{r}{\Cosp(\lim)} & \Cosp(\Pro(\catsingle{S}))
        \end{tikzcd}
    \end{center}
\end{construction}
We have fairly understood how the functors $\inc$ and $\lim$ in the above behave on objects. We now explain what $s_!$ does on objects. Fix $X_{\bullet,\bullet}\in \Fun(\mathbb{I},\catsingle{S})$. Then, by the pointwise formula for left Kan extension, we obtain
$$s_! X(n)\simeq \colim_{[i,j],\, i\to n}X_{i,j}\simeq \colim_{[i,j],\,i\geq n}X_{i,j}$$
Consider the inclusion $\iota\colon\{[n,m]\} \to \{[i,j], i\geq n\}$; this is clearly cofinal. Indeed, for a fixed $[i,j]$, the slice $[i,j]_{/\iota}$ is the set $\{[n,m],\,n\leq i\leq j\leq m\}$ which clearly has an initial object, hence contractible. Thus,
$$s_! X(n)\simeq \colim_{m\geq n}X_{n,m}$$
where the maps in the diagram are given by
$$X_{n,n} \xrightarrow[]{X(\mathrm{id},n\to n+1)} X_{n,n+1} \xrightarrow[]{X(\id,n+1\to n+2)}X_{n,n+2} \to \cdots$$
We therefore suggestively denote $s_! X(n)$ by $X_{n,\infty}$

In order to decode the rather confusing definition of the above functor, we illustrate what it does on objects through the following diagram.
\begin{center}
    \begin{tikzcd}[column sep={3.6em,between origins}, row sep={3.3em,between origins},
  po/.style={phantom, very near start}]
& & & & |[blue]| X_{0,\infty}
& & |[blue]| X_{1,\infty} \arrow[ll, blue]
& & |[blue]| X_{2,\infty} \arrow[ll, blue]
& & |[blue]| \cdots \arrow[ll, blue] \\
& & & X_{0,3} \arrow[ur, dashed, blue]
  \arrow[dd, po, "\ulcorner"{rotate=135}]
& \vdots
& X_{1,4} \arrow[ur, dashed, blue]
  \arrow[dd, po, "\ulcorner"{rotate=135}]
& & X_{2,5} \arrow[ur, dashed, blue]
  \arrow[dd, po, "\ulcorner"{rotate=135}]
& & \cdots & \\
& & X_{0,2} \arrow[ur]
  \arrow[dd, po, "\ulcorner"{rotate=135}]
& & X_{1,3} \arrow[ul] \arrow[ur]
  \arrow[dd, po, "\ulcorner"{rotate=135}]
& & X_{2,4} \arrow[ul] \arrow[ur]
  \arrow[dd, po, "\ulcorner"{rotate=135}]
& & & & \\
& |[red]| X_{0,1} \arrow[ur]
& & |[red]| X_{1,2} \arrow[ul] \arrow[ur]
& & |[red]| X_{2,3} \arrow[ul] \arrow[ur]
& & |[red]| X_{3,4} \arrow[ul]
& & |[red]| \cdots & \\
|[red]| X_0 \arrow[ur, red]
& & |[red]| X_1 \arrow[ul, red] \arrow[ur, red]
& & |[red]| X_2 \arrow[ul, red] \arrow[ur, red]
& & |[red]| X_3 \arrow[ul, red] \arrow[ur, red]
& & |[red]| X_4 \arrow[ul, red]
& &
\end{tikzcd}
\end{center}
The red part of the above diagram is the same as the datum of a functor $\mathbb{N} \to \Cosp(\catsingle{S})$; namely, a sequence of spaces $X_i$, and cospans $X_{i,i+1}$ between them. The black part of the diagram is taking iterated pushouts. The blue part is then the diagonal colimits, and we observe that the inverse system $(X_{0,\infty}\leftarrow X_{1,\infty} \leftarrow X_{2,\infty}\leftarrow\cdots)$ in blue is the object in $\Pro(\catsingle{S})$ that the functor $X_\bullet\colon \mathbb{N} \to \Cosp(\catsingle{S})$ is mapped to by the tail functor.

The following is a crucial lemma, allowing us to ensure that the tail functor descends to a functor on $\PSh^\mathbb{N}(\Cosp(\catsingle{S}))$:
\begin{lem}\label{lem:tail-functor}
    The functor $\mathrm{tail}\colon \Fun(\mathbb{N},\Cosp(\catsingle{S})) \to \Cosp(\Pro(\catsingle{S}))$ sends a cofinal re-indexing morphism, i.e. a morphism in $W$ as in Theorem \ref{thm:localisation}, to an equivalence in $\Cosp(\Pro(\catsingle{S}))$. As a consequence of Theorem \ref{thm:localisation}, $\mathrm{tail}$ descends to a functor
    $$\tau\colon \PSh^\mathbb{N}(\Cosp(\catsingle{S})) \to \Cosp(\Pro(\catsingle{S}))$$
\end{lem}
\begin{proof}
    We fix a morphism in $W$, which by definition is the following: fix a strictly increasing function $f \colon \mathbb{N} \to \mathbb{N}$; then, for $X \in \Fun(\mathbb{N},\Cosp(\catsingle{S}))$, we consider the functor $f^\ast X \coloneqq X \circ f$. The morphism we fix in $W$ is then the morphism $\omega \colon X \to f^\ast X$. We claim that $\omega$ is sent to an equivalence in $\Cosp(\Pro(\catsingle{S}))$ via the tail functor. Consider the cospan between $X_i$ and $X_{f(i)}$, namely $X_i \rightarrow X_{i,f(i)} \leftarrow X_{f(i)}$. We then obtain a commutative diagram
    \begin{center}
        \begin{tikzcd}
            X_i \arrow{r} \arrow{d} & X_{i,f(i)} \arrow{d} & \arrow{l} X_{f(i)} \arrow{d}
            \\
            X_{i,i+1} \arrow{r}\arrow{d} & X_{i,f(i+1)} \arrow{d} & X_{f(i),f(i+1)} \arrow{l} \arrow{d}
            \\
            \vdots & \vdots & \vdots
        \end{tikzcd}
    \end{center}
    and the induced cospan of $\Pro(\catsingle{S})$ is induced by applying vertical colimits, namely
    $$\{X_{i,\infty}\}_i \rightarrow \{\colim_{n}X_{i,f(n)}\}_i \leftarrow \{\colim_n X_{f(i),f(n)}\}_i$$
   The left hand map is an equivalence objectwise (by cofinality of $f$), and hence an equivalence of the associated $\Pro$-systems. The right hand side is the reindexing of a tower along a cofinal map, and hence induces an equivalence of $\Pro$-objects. Thus we obtain a cospan of $\Pro(\catsingle{S})$, where both legs of the cospans are equivalences; these are the equivalences in $\Cosp(\Pro(\catsingle{S}))$. That is, $\omega$ is mapped to an equivalence in $\Cosp(\Pro(\catsingle{S}))$, as claimed.
\end{proof}

\subsection{The $\proEnd$ functor}
We are now in a position to define the required functor recording the $\Pro$-end of a manifold, namely $\proEnd\colon \ManInf_d \to \Cosp(\Pro(\catsingle{S}))$. A starting step is the following
\begin{lem}\label{lem:bordism-to-cosp}
    For all $d\in \mathbb{N}$, there is a functor $u \colon \BordInf_d \to \Cosp(\catsingle{S})$, sending a closed $(d-1)$-dimensional manifold to the homotopy type, and a bordism $W\colon M\rightsquigarrow N$ to the cospan in spaces $(M \rightarrow W\leftarrow N)$, where the legs are the underlying maps of the boundary embeddings $M\sqcup N \cong \partial W \hookrightarrow W$.
\end{lem}
\begin{proof}
    We consider the functor $\fr(-)\colon \ncBordInf_d \to \COSPAN^+(\catsingle{S}_{/\BO(d)})$ defined in \cite[p. 46]{discstructurespace}. Here, $\ncBordInf_d$ is the double $\infty$-bordism category of non-compact bordisms constructed in \emph{loc. cit.}, and where $\COSPAN^+(\catsingle{S}_{/\BO(d)})$ refers to the double $\infty$-category of cospans as in \emph{loc. cit.}, §2.10, based on \cite{HaugsengIteratedSpans}. Furthermore, in \cite{discstructurespace}, there is a functor $\BordInf_d \to \ncBordInf_d$, where now $\BordInf_d$ is viewed as a double category in the obvious manner: an $(\infty,1)$-category is a category object in $\catsingle{S}$, namely a simplicial object in $\catsingle{S}$ satisfying the Segal and completeness conditions; a double category is a category object in $\CatInf$, and there is a fully faithful functor $\catsingle{S} \hookrightarrow \CatInf$, which, applied levelwise to $\Fun(\Delta^\op,\catsingle{S})$, yields a double category. Restriction of $\fr(-)\colon \ncBordInf_d \to \COSPAN^+(\catsingle{S}_{/\BO(d)})$ along the subcategory inclusion $\BordInf_d \to \ncBordInf_d$ thus gives rise to a functor $\BordInf_d \to \mathrm{cospan}^+_1(\catsingle{S}_{/\BO(d)})$, which by abuse of notation we also denote $\fr(-)$: here, $\mathrm{cospan}_1^+$ denotes the $(\infty,1)$ category of cospans as follows from \cite[Definition 5.8]{HaugsengIteratedSpans}. It remains to show that $\mathrm{cospan}_1^+$ above, and $\Cosp$ as in \cite{Spans}, are equivalent categories. To see this, we begin by noting that they are both defined as certain simplicial objects in $\catsingle{S}$ that satisfy the Segal and completeness conditions. We claim that the simplicial objects themselves are naturally equivalent. Indeed, if $\catsingle{C}$ is a category admitting pullbacks, the definition of $\mathrm{span}(\catsingle{C})$ in \cite{Spans} spells out to
    $$\mathrm{Span}(\catsingle{C})_{[n]}\coloneqq \Map_{\mathrm{AdTrip}}(\tw^r([n]),\catsingle{C})$$
    where $\catsingle{C}$ is given its adequate triple structure via $(\catsingle{C},\catsingle{C},\catsingle{C})$, and $\tw^r([n])$ its adequate triple structure as in \cite[Example 2.9]{Spans}. The condition on sending ambigressive pullbacks to pullbacks thus spells out to the following: functors $F\colon \tw^r([n])\to \catsingle{C}$, such that the natural map $F(i,\ell) \to F(i,j)\times_{F(k,j)}F(k,\ell)$ is an equivalence, for $i\leq k\leq j\leq \ell$. By \cite[Remark 5.7]{HaugsengIteratedSpans}, this description of the underlying simplicial space of $\mathrm{Span}(\catsingle{C})$ is equivalent to the one in \cite[Definition 5.8]{HaugsengIteratedSpans}.

    We are now equipped with a functor $\fr(-)\colon\BordInf_d \to \Cosp(\catsingle{S}_{/\BO(d)})$. As the forgetful functor $\catsingle{S}_{/\BO(d)}\to \catsingle{S}$ preserves colimits, it induces a functor $\Cosp(\catsingle{S}_{/\BO(d)}) \to \Cosp(\catsingle{S})$. Precomposing with $\fr$, we obtain our desired functor $u\colon \BordInf_d \to \Cosp(\catsingle{S})$.
\end{proof}
We thus define $\proEnd$ as follows:
\begin{dfn}
    Define $\proEnd \colon \ManInf_d \to \Cosp(\Pro(\catsingle{S}))$ as the following composite:
    $$\ManInf_d \xrightarrow[]{\bordEnd}\PSh^\mathbb{N}(\BordInf_d) \xrightarrow[]{\PSh^\mathbb{N}(u)}\PSh^\mathbb{N}(\Cosp(\catsingle{S})) \xrightarrow[]{\tau}\Cosp(\Pro(\catsingle{S}))$$
    where
    \begin{itemize}
        \item $\bordEnd$ is the functor of Definition \ref{def:boundary-bordismends-functors};
        \item $\PSh^\mathbb{N}(u)$ is the functor obtained by applying $\PSh^\mathbb{N}$ to the functor $u$ of Lemma \ref{lem:bordism-to-cosp};
        \item $\tau$ is the functor of Lemma \ref{lem:tail-functor}.
    \end{itemize}
\end{dfn}
In particular, from the mere definition, we have achieved the following:
\begin{thm}\label{thm:bordEnd-proEnd}
    There is a commutative triangle in $\CatInf$
    \begin{center}
        \begin{tikzcd}
            & \PSh^\mathbb{N}(\BordInf_d) \arrow{dd}{\tau\circ\PSh^\mathbb{N}(u)}
            \\
            \ManInf_d \arrow{ur}{\bordEnd} \arrow{dr}[swap]{\proEnd}
            \\
            & \Cosp(\Pro(\catsingle{S}))
        \end{tikzcd}
    \end{center}
    Furthermore, for $W\in\OpenMan_d$, and a choice of a compact exhaustion $\{K_i\}$ thereof, we have the equivalences
    \begin{itemize}
        \item $\bordEnd(W)\simeq \{\partial K_i\}$ in $ \PSh^\mathbb{N}(\BordInf_d)$;
        \item $\proEnd(W) \simeq \{W\setminus K_i\}$ in $\Cosp(\Pro(\catsingle{S}))$
    \end{itemize}
\end{thm}

\begin{rem}
    For a category $\catsingle{C}$ with pushouts, there is a map $\catsingle{C} \to \Cosp(\catsingle{C})$, sending an object to itself, and a morphism $f\colon c \to d$ to the cospan $c\xrightarrow[]{f}d \xleftarrow[]{\id}d$. On groupoid cores, the above functor is an equivalence (for instance, see the discussion following Theorem 2.13 in \cite{Spans}). Thus, applying groupoid cores to the above theorem, we obtain a commutative diagram
    \begin{center}
        \begin{tikzcd}
            & \PSh^\mathbb{N}(\BordInf_d)^\simeq \arrow{dd}{\tau\circ\PSh^\mathbb{N}(u)}
            \\
            \left(\ManInf_d\right)^\simeq \arrow{ur}{\bordEnd} \arrow{dr}[swap]{\proEnd}
            \\
            & \Pro(\catsingle{S})^\simeq
        \end{tikzcd}
    \end{center}
\end{rem}

\bibliographystyle{amsalpha}
\bibliography{literature}

@misc{AdelaJan,
      title={Open 2D TFTs admit initial open-closed extensions}, 
      author={Shaul Barkan and Jan Steinebrunner and Adela YiYu Zhang},
      year={2025},
      eprint={2509.02553},
      archivePrefix={arXiv},
      primaryClass={math.AT},
      url={https://arxiv.org/abs/2509.02553}, 
      howpublished={arxiv:2509.02553}
}

@book {SiebenmannThesis,
    AUTHOR = {Siebenmann, Laurence Carl},
     TITLE = {The obstruction to finding a boundary for an
              open manifold of dimension greater than five},
      NOTE = {Thesis (Ph.D.)--Princeton University},
 PUBLISHER = {ProQuest LLC, Ann Arbor, MI},
      YEAR = {1965},
     PAGES = {152},
   MRCLASS = {99-05},
  MRNUMBER = {2615648},
       URL = {http://gateway.proquest.com/openurl?url_ver=Z39.88-2004&rft_val_fmt=info:ofi/fmt:kev:mtx:dissertation&res_dat=xri:pqdiss&rft_dat=xri:pqdiss:6605012},
}

@article {Weiss,
    AUTHOR = {Weiss, Michael},
     TITLE = {Embeddings from the point of view of immersion theory. {I}},
   JOURNAL = {Geom. Topol.},
  FJOURNAL = {Geometry and Topology},
    VOLUME = {3},
      YEAR = {1999},
     PAGES = {67--101},
      ISSN = {1465-3060,1364-0380},
   MRCLASS = {57R40 (57R42)},
  MRNUMBER = {1694812},
MRREVIEWER = {R.\ E.\ Stong},
       DOI = {10.2140/gt.1999.3.67},
       URL = {https://doi.org/10.2140/gt.1999.3.67},
}

@book {HTT,
    AUTHOR = {Lurie, Jacob},
     TITLE = {Higher topos theory},
    SERIES = {Annals of Mathematics Studies},
    VOLUME = {170},
 PUBLISHER = {Princeton University Press, Princeton, NJ},
      YEAR = {2009},
     PAGES = {xviii+925},
      ISBN = {978-0-691-14049-0; 0-691-14049-9},
   MRCLASS = {18-02 (18B25 18E35 18G30 18G55 55U40)},
  MRNUMBER = {2522659},
MRREVIEWER = {Mark\ Hovey},
       DOI = {10.1515/9781400830558},
       URL = {https://doi.org/10.1515/9781400830558},
}

@article {Rezk,
    AUTHOR = {Rezk, Charles},
     TITLE = {Free colimit completion in {$\infty$}-categories},
   JOURNAL = {Homology Homotopy Appl.},
  FJOURNAL = {Homology, Homotopy and Applications},
    VOLUME = {27},
      YEAR = {2025},
    NUMBER = {1},
     PAGES = {275--292},
      ISSN = {1532-0073,1532-0081},
   MRCLASS = {18N60 (18A30 55U35)},
  MRNUMBER = {4889274},
MRREVIEWER = {Martin\ Frankland},
       DOI = {10.4310/hha.2025.v27.n1.a15},
       URL = {https://doi.org/10.4310/hha.2025.v27.n1.a15},
}

@misc{derivedmappingspacesinftycategories,
      title={Derived mapping spaces of $\infty$-categories}, 
      author={Kensuke Arakawa and Daniel Carranza and Chris Kapulkin},
      year={2025},
      eprint={2509.10288},
      archivePrefix={arXiv},
      primaryClass={math.AT},
      url={https://arxiv.org/abs/2509.10288},
      howpublished={arxiv:2509.10288}
}

@misc{efimov,
      title={Localizing invariants of inverse limits}, 
      author={Alexander I. Efimov},
      year={2025},
      eprint={2502.04123},
      archivePrefix={arXiv},
      primaryClass={math.KT},
      url={https://arxiv.org/abs/2502.04123},
      howpublished={arxiv:2502.04123}
}

@article {GepnerHaugsengNikolaus,
    AUTHOR = {Gepner, David and Haugseng, Rune and Nikolaus, Thomas},
     TITLE = {Lax colimits and free fibrations in {$\infty$}-categories},
   JOURNAL = {Doc. Math.},
  FJOURNAL = {Documenta Mathematica},
    VOLUME = {22},
      YEAR = {2017},
     PAGES = {1225--1266},
      ISSN = {1431-0635,1431-0643},
   MRCLASS = {18D30 (18A30)},
  MRNUMBER = {3690268},
MRREVIEWER = {Josep\ Elgueta},
}

@book {Hirsch,
    AUTHOR = {Hirsch, Morris W.},
     TITLE = {Differential topology},
    SERIES = {Graduate Texts in Mathematics},
    VOLUME = {33},
      NOTE = {Corrected reprint of the 1976 original},
 PUBLISHER = {Springer-Verlag, New York},
      YEAR = {1994},
     PAGES = {x+222},
      ISBN = {0-387-90148-5},
   MRCLASS = {57-01 (58-01)},
  MRNUMBER = {1336822},
}

@article {Spans,
    AUTHOR = {Haugseng, Rune and Hebestreit, Fabian and Linskens, Sil and
              Nuiten, Joost},
     TITLE = {Two-variable fibrations, factorisation systems and
              {$\infty$}-categories of spans},
   JOURNAL = {Forum Math. Sigma},
  FJOURNAL = {Forum of Mathematics. Sigma},
    VOLUME = {11},
      YEAR = {2023},
     PAGES = {Paper No. e111, 70},
      ISSN = {2050-5094},
   MRCLASS = {18N70 (18N60 18N65)},
  MRNUMBER = {4676218},
MRREVIEWER = {Fredrik\ Nordvall Forsberg},
       DOI = {10.1017/fms.2023.107},
       URL = {https://doi.org/10.1017/fms.2023.107},
}

@book {Cardim,
    AUTHOR = {Cardim, Nancy de Souza},
     TITLE = {Embeddings of open manifolds},
      NOTE = {Thesis (Ph.D.)--University of Notre Dame},
 PUBLISHER = {ProQuest LLC, Ann Arbor, MI},
      YEAR = {1993},
     PAGES = {54},
   MRCLASS = {99-05},
  MRNUMBER = {2689783},
       URL ={http://gateway.proquest.com/openurl?url_ver=Z39.88-2004&rft_val_fmt=info:ofi/fmt:kev:mtx:dissertation&res_dat=xri:pqdiss&rft_dat=xri:pqdiss:9400152},
}

@misc{Agol,
  title = {Isotopy inverse embeddings vs. diffeomorphisms},
  author = {Ian Agol},
  howpublished = {MathOverflow answer},
  note = {URL: https://mathoverflow.net/questions/100498/isotopy-inverse-embeddings-vs-diffeomorphisms},
  eprint = {https://mathoverflow.net},
  url = {https://mathoverflow.net}
}

@book {MilnorMorse,
    AUTHOR = {Milnor, J.},
     TITLE = {Morse theory},
    SERIES = {Annals of Mathematics Studies},
    VOLUME = {No. 51},
      NOTE = {Based on lecture notes by M. Spivak and R. Wells},
 PUBLISHER = {Princeton University Press, Princeton, NJ},
      YEAR = {1963},
     PAGES = {vi+153},
   MRCLASS = {57.50 (53.72)},
  MRNUMBER = {163331},
MRREVIEWER = {H.\ I.\ Levine},
}

@article {discstructurespace,
    AUTHOR = {Krannich, Manuel and Kupers, Alexander},
     TITLE = {The disc-structure space},
   JOURNAL = {Forum Math. Pi},
  FJOURNAL = {Forum of Mathematics. Pi},
    VOLUME = {12},
      YEAR = {2024},
     PAGES = {Paper No. e26, 98},
      ISSN = {2050-5086},
   MRCLASS = {57S05 (18F50 18N60 18N70 58D10)},
  MRNUMBER = {4838966},
MRREVIEWER = {Marja\ K.\ Kankaanrinta},
       DOI = {10.1017/fmp.2024.25},
       URL = {https://doi.org/10.1017/fmp.2024.25},
}

@incollection {Stallings,
    AUTHOR = {Stallings, John R.},
     TITLE = {On infinite processes leading to differentiability in the
              complement of a point},
 BOOKTITLE = {Differential and {C}ombinatorial {T}opology ({A} {S}ymposium
              in {H}onor of {M}arston {M}orse)},
     PAGES = {245--254},
 PUBLISHER = {Princeton Univ. Press, Princeton, NJ},
      YEAR = {1965},
   MRCLASS = {57.10},
  MRNUMBER = {180983},
MRREVIEWER = {N.\ Kuiper},
}

@article {ManuelSanderInfinityOperadic,
    AUTHOR = {Krannich, Manuel and Kupers, Alexander},
     TITLE = {Infinity-operadic foundations for embedding calculus},
   JOURNAL = {J. Topol.},
  FJOURNAL = {Journal of Topology},
    VOLUME = {19},
      YEAR = {2026},
    NUMBER = {2},
     PAGES = {Paper No. e70071, 136},
      ISSN = {1753-8416,1753-8424},
   MRCLASS = {18F50 (18N70 57S05 58D10)},
  MRNUMBER = {5062094},
       DOI = {10.1112/topo.70071},
       URL = {https://doi.org/10.1112/topo.70071},
}

@article {HaugsengIteratedSpans,
    AUTHOR = {Haugseng, Rune},
     TITLE = {Iterated spans and classical topological field theories},
   JOURNAL = {Math. Z.},
  FJOURNAL = {Mathematische Zeitschrift},
    VOLUME = {289},
      YEAR = {2018},
    NUMBER = {3-4},
     PAGES = {1427--1488},
      ISSN = {0025-5874,1432-1823},
   MRCLASS = {18D05 (14A20 14F05 81R10)},
  MRNUMBER = {3830256},
MRREVIEWER = {Benjamin\ Cooper},
       DOI = {10.1007/s00209-017-2005-x},
       URL = {https://doi.org/10.1007/s00209-017-2005-x},
}

@article {diffcompactsupport,
    AUTHOR = {Banakh, Taras and Yagasaki, Tatsuhiko},
     TITLE = {Diffeomorphism groups of non-compact manifolds endowed with
              the {W}hitney {$C^\infty$}-topology},
   JOURNAL = {Topology Appl.},
  FJOURNAL = {Topology and its Applications},
    VOLUME = {179},
      YEAR = {2015},
     PAGES = {51--61},
      ISSN = {0166-8641,1879-3207},
   MRCLASS = {57N20 (46A13 57S05 58D05)},
  MRNUMBER = {3270925},
MRREVIEWER = {Daniel\ Belti\c t\u a},
       DOI = {10.1016/j.topol.2014.08.015},
       URL = {https://doi.org/10.1016/j.topol.2014.08.015},
}

@book {ANR,
    AUTHOR = {Aguilar, Marcelo and Gitler, Samuel and Prieto, Carlos},
     TITLE = {Algebraic topology from a homotopical viewpoint},
    SERIES = {Universitext},
      NOTE = {Translated from the Spanish by Stephen Bruce Sontz},
 PUBLISHER = {Springer-Verlag, New York},
      YEAR = {2002},
     PAGES = {xxx+478},
      ISBN = {0-387-95450-3},
   MRCLASS = {55-01},
  MRNUMBER = {1908260},
MRREVIEWER = {Donald\ M.\ Davis},
       DOI = {10.1007/b97586},
       URL = {https://doi.org/10.1007/b97586},
}

@article {Palais,
    AUTHOR = {Palais, Richard S.},
     TITLE = {Homotopy theory of infinite dimensional manifolds},
   JOURNAL = {Topology},
  FJOURNAL = {Topology. An International Journal of Mathematics},
    VOLUME = {5},
      YEAR = {1966},
     PAGES = {1--16},
      ISSN = {0040-9383},
   MRCLASS = {57.55 (55.40)},
  MRNUMBER = {189028},
MRREVIEWER = {A.\ Dold},
       DOI = {10.1016/0040-9383(66)90002-4},
       URL = {https://doi.org/10.1016/0040-9383(66)90002-4},
}

\end{document}